\documentclass[11pt,twoside,english]{myamsart}
\advance\oddsidemargin by -1.0cm
\advance\evensidemargin by -1.0cm
\advance\topmargin by -1.0cm 
\usepackage{amssymb}
\usepackage{babel}
\usepackage{amsmath}
\usepackage{amscd}   
\usepackage{epsfig}  
\usepackage{rotating}
\let\mathg\mathfrak
\theoremstyle{plain}
\newtheorem{cor}{Corollary}[section]
\newtheorem{lem}{Lemma}[section]
\newtheorem{thm}{Theorem}[section]
\newtheorem{prop}{Proposition}[section]
\theoremstyle{definition}
\newtheorem{exa}{Example}[section]
\newtheorem{NB}{Remark}[section]
\newtheorem{dfn}{Definition}[section]

\renewcommand{\labelenumi}{$(\theenumi)$}
\newcommand{\bdm}{\begin{displaymath}}
\newcommand{\edm}{\end{displaymath}}
\newcommand{\be}{\begin{equation}}
\newcommand{\ee}{\end{equation}}
\newcommand{\ba}[1]{\begin{array}{#1}}
\newcommand{\ea}{\end{array}}

\newcommand{\btab}{\begin{tabular}}
\newcommand{\etab}{\end{tabular}}

\newcommand{\x}{\times}
\newcommand{\op}{\oplus}
\newcommand{\ox}{\otimes}

\newcommand{\ra}{\rightarrow}
\newcommand{\lra}{\longrightarrow}

\newcommand{\tr}{\ensuremath{\mathrm{tr}}}

\newcommand{\R}{\ensuremath{\mathbb{R}}}

\newcommand{\Z}{\ensuremath{\mathbb{Z}}}

\newcommand{\D}{\ensuremath{\mathcal{D}}} 

\newcommand{\T}{\ensuremath{\mathrm{T}}}

\newcommand{\eps}{\ensuremath{\varepsilon}} 
    
\newcommand{\vrho}{\ensuremath{\varrho}}
\newcommand{\hut}{\wedge}

\newcommand{\Scal}{\ensuremath{\mathrm{Scal}}}

\newcommand{\Ad}{\ensuremath{\mathrm{Ad}\,}}

\newcommand{\diag}{\ensuremath{\mathrm{diag}}}

\newcommand{\GL}{\ensuremath{\mathrm{GL}}}

\newcommand{\su}{\ensuremath{\mathg{su}}}
\newcommand{\SU}{\ensuremath{\mathrm{SU}}}

\newcommand{\so}{\ensuremath{\mathg{so}}}
\newcommand{\SO}{\ensuremath{\mathrm{SO}}}
\newcommand{\Spin}{\ensuremath{\mathrm{Spin}}}
\newcommand{\spin}{\ensuremath{\mathg{spin}}}
\newcommand{\Orth}{\ensuremath{\mathrm{O}}}

\newcommand{\g}{\ensuremath{\mathfrak{g}}}
\newcommand{\h}{\ensuremath{\mathfrak{h}}}
\newcommand{\m}{\ensuremath{\mathfrak{m}}}

\begin{document}
\def\haken{\mathbin{\hbox to 6pt{%
                 \vrule height0.4pt width5pt depth0pt
                 \kern-.4pt
                 \vrule height6pt width0.4pt depth0pt\hss}}}
    \let \hook\intprod
\setcounter{equation}{0}
%
%
\thispagestyle{empty}
%
\date{\today}
\title{On the holonomy of connections with skew-symmetric torsion}
%
%
%
\author{Ilka Agricola}
\author{Thomas Friedrich}
\address{\hspace{-5mm} 
{\normalfont\ttfamily agricola@mathematik.hu-berlin.de}\newline
{\normalfont\ttfamily friedric@mathematik.hu-berlin.de}\newline
Institut f\"ur Mathematik \newline
Humboldt-Universit\"at zu Berlin\newline
Unter den Linden 6\newline
Sitz: John-von-Neumann-Haus, Adlershof\newline
D-10099 Berlin, Germany}
%
\thanks{Supported by the SFB 288 "Differential geometry
and quantum physics" of the DFG}
\subjclass[2000]{Primary 53 C 25; Secondary 81 T 30}
\begin{abstract}
We investigate the holonomy group of a linear metric connection with
skew-sym\-me\-tric torsion. In case of the euclidian space and a constant torsion
form this group is always semisimple.
It does not preserve any non-degenerated $2$-form or any spinor.
Suitable integral formulas allow us to prove similar properties 
in case of a compact
Riemannian manifold equipped with a metric connection of skew-symmetric
torsion. On the Aloff-Wallach space $N(1,1)$ we construct families
of connections admitting parallel spinors. Furthermore, we investigate
the geometry of these connections as well as the geometry of the
underlying Riemannian metric. Finally, we prove that any $7$-dimensional
$3$-Sasakian manifold admits $\mathbb{P}^2$-parameter families of
linear metric connections and spinorial connections defined by $4$-forms
with  parallel spinors.
\end{abstract}
\maketitle
\tableofcontents
\pagestyle{headings}
%
%
%
\section{Introduction}\noindent
This paper treats the geometry of metric invariant connections with 
skew-sym\-me\-tric torsion, as they became recently of interest in string 
theory and special geometries.  The notion of torsion of a connection 
was invented by Elie Cartan,
and appeared for the first time in a short note at the Acad\'emie des Sciences
de Paris in 1922 (see \cite{Cartan1})\footnote{We thank Andrzej Trautman 
for drawing our attention to these papers by Cartan -- see \cite{Tr}.}.
Although it contains no formulas, 
Cartan observes that such a connection may or may not preserve geodesics,
and turns his attention first to those who actually do so. In this sense,
E.~Cartan was the first to investigate this class of 
connections. At that time, it was not yet customary -- as it became later
in the second half of the 20th century -- to assign to a Riemannian manifold
only its Levi-Civita connection. Rather, Cartan demands (see \cite{Cartan4}):
\begin{quote}
\emph{\'Etant donn\'e une vari\'et\'e plong\'ee dans l'espace affine (ou
projectif, ou conforme etc.), attribuer \`a cette vari\'et\'e la
connexion affine (ou projective, ou conforme etc.) qui rende le plus
simplement compte des relations de cette vari\'et\'e avec l'espace ambiant.} 
\end{quote}
He then goes on to explain in very general terms how the connection should
be \emph{adapted} to the geometry under consideration. This point of
view should be taken into account in Riemannian geometry, too. The
canonical connection of a naturally reductive Riemannian space is a first
example (see \cite{Agri}). Moreover, we know many non integrable 
geometric structures 
on Riemannian manifolds admitting a unique metric connection preserving the 
structure and with non vanishing
skew-symmetric torsion (see \cite{Friedrich&I1}, \cite{Fri2}). Following 
Cartan
as well as the idea that torsion forms  are candidates
for the so called $B$-field in string theory, the geometry of these connections
deserves systematic investigation. Basically, there are no general results 
concerning the
holonomy group of connections with torsion. The question whether or
not a connection of that type admits parallel tensor fields 
differs radically from the corresponding problem  for the Levi-Civita 
connection. In particular, one is interested in  the existence of 
parallel spinor fields,
since they are interpreted in string theory as  supersymmetries of the
model.

\smallskip\noindent
The paper is organized as follows. In  Section~$2$,  we discuss once again 
some basic results motivating the role of metric connections with 
skew-symmetric torsion.
In Sections~$3$  and $4$,  we study the linear case, i.e.,
euclidian space equipped with a constant torsion form $\T$. The holonomy
algebra $\g_{\T}^*$ of the corresponding linear connection has some
remarkable properties. For any $3$-form, $\g_{\T}^*$ is a semisimple
Lie algebra. Moreover, it cannot preserve a non degenerate $2$-form
or a spinor. On the other side, many representations of a compact,
semisimple Lie algebra occur as the holonomy algebra of some $3$-form,
for example the adjoint representation can be realized in this way. We
introduce an  obstruction for a Lie algebra representation 
to be the holonomy algebra of some $3$-form and show on an example
how it may be used to rule  out some representations. In particular, 
the unique, irreducible $16$-dimensional representation of the 
algebra $\spin(9)$ cannot be the holonomy algebra of some $3$-form. 
Forms of higher degree
than three do not occur for linear connections, but they define spinorial
connections. In the euclidian case we introduce their spinorial
holonomy algebra as a Lie subalgebra of the Clifford algebra. In all
examples discussed, this algebra turns out to be perfect.

\smallskip\noindent
In Section $5$ and $6$, we generalize the algebraic results to the case
of a Riemannian manifold $(M^n, g , \T)$ with a metric connection $\nabla$.
In particular, we are interested in the question whether or not the
$\nabla$-holonomy group preserves a spinor field. In the compact case,
we prove that if the scalar curvature $\mathrm{Scal}^g \leq 0$ is non
positive and if the torsion form is closed, $d \T = 0$,  any
$\nabla$-parallel spinor is Riemannian parallel and $\T = 0$ vanishes.
Here we use an integral formula for the square of the Dirac operator depending
on the connection. The main point is that the formula becomes
simple if one compares the Dirac operator corresponding the connection
with torsion form $\T$ with the spinorial Laplace operator corresponding
to the connection with torsion form $3\cdot\T$. This effect has been
observed in the literature at several places, in particular by
Bismut (see \cite{Bis}) and, in the homogeneous case, by 
Agricola (see \cite{Agri}). We explore the corresponding integral
formula and study the space of parallel spinors.

\smallskip\noindent
In Section $7$ we discuss, for a given triple $(M^n,g,\T)$, 
the whole family $\nabla^s$ of metric connections with torsion form
$s \cdot \T$. In the generic case, the existence of a $\nabla^s$-parallel 
spinor restricts the possible parameter $s$ via
a polynomial equation. Consequently, in the generic case, at most a finite
number of connections in the family admits parallel spinors. Some
simple examples show that sometimes two connections 
really admit parallel spinors. Moreover,
our integral formulas prove that, on a compact manifold, basically 
only three parameters are possible. In case that the torsion form
is associated with a special non integrable geometry, 
the connection $\nabla^s$ with a parallel
spinor is sometimes unique. A result of that type 
requires additional informations concerning the underlying geometry. We prove
it for $5$-dimensional Sasakian manifolds equipped with their 
canonical connection. 

\smallskip\noindent
In Section $8$ we construct, on the Aloff-Wallach manifold 
$N(1,1) =\SU(3)/S^1$,  a two-parameter family of metrics 
that admits two inequivalent cocalibrated $G_2$-structures. 
Moreover, we investigate the torsion forms of their unique connections 
as well as other geometric data of these connections. Our approach is different from the usual one (see \cite{CMS}). First we construct $3$-forms with
parallel spinors on $N(1,1)$. The underlying $\mathrm{G}_2$-structure
is cocalibrated and many of the geometric data are encoded into
the torsion $3$-form we started with. Moreover, we are interested not only
in the type of the $\mathrm{G}_2$-structure, but mainly in the geometry
of the unique connection preserving this structure. The same method is then
applied in order to construct spinorial connections defined by
$4$-form and admitting parallel spinor fields. Some of these connections
are closely related to the $3$-Sasakian structure of $N(1,1)$. In section
$9$ we generalize these examples. Indeed, we are able to construct,
for any $7$-dimensional $3$-Sasakian manifold, a canonical 
$\mathbb{P}^2$-parameter family of $3$- and $4$-forms such that the 
underlying linear or spinorial connection  admits parallel
spinors.
 
%
%
\section{The eight classes of linear connections with torsion}
\label{classes}\noindent
We begin by an
elementary, yet enlightening investigation of geometric torsion tensors.  
Consider a Riemannian manifold $(M^n, g)$. In a point,
the difference between its Levi-Civita connection $\nabla^g$ and any linear 
connection $\nabla$ is  a $(2,1)$-tensor $A$,
\bdm 
\nabla_X Y\ =\ \nabla^g_X Y + A(X,Y),\quad X,Y \in TM\,.
\edm
The vanishing of the symmetric or the antisymmetric part of $A$ 
has an immediate geometric interpretation:
\begin{lem}\label{same-geodesics}
The connection $\nabla$ is torsion-free if and only if $A$ is symmetric.
The connection $\nabla$ has the same geodesics as the Levi-Civita connection 
$\nabla^g$ if and only if $A$ is antisymmetric.
\end{lem}
\begin{proof}
The torsion $\T$ of $\nabla$ is 
\bdm
\T(X,Y)\ :=\ \nabla_X Y - \nabla _Y X -[X,Y]\ =\ A(X,Y) - A(Y,X) \, ,
\edm
since $\nabla^g$ is torsion-free. Hence the first claim follows.
For the second, consider a curve $\gamma$ through a point $p$, and
set $X:=\dot{\gamma}(p)$. Then
\bdm
\nabla_{X} X\ = \  
\nabla^g_{X}X+ A(X,X)\, ,
\edm
and hence $\nabla_{X} X$ coincides with $\nabla_{X}^g X$ if and only if
$A$ is skew-symmetric.
\end{proof}
\noindent
Following Cartan (see \cite[p.51]{Cartan3}), we study the algebraic 
types of the torsion tensor for a metric connection.
Denote by the same symbol the $(3,0)$-tensors
derived from $A,\T$ by contraction with the metric,
\bdm
A(X,Y,Z)\ :=\ g( A(X,Y),Z)\, , \quad \T(X,Y,Z) \ := \ g(\T(X,Y),Z) \, .
\edm
We identify $TM$ with $TM^*$ via the metric from now on. 
Let $\mathcal{T}$ be the $n^2(n-1)/2$-dimensional space of all 
possible torsion tensors,
\bdm
\mathcal{T}\ =\ \{\T\in\ox^3 TM \ | \ \T(X,Y,Z)= - \, \T(Y,X,Z) \}
\ \cong \ \wedge^2 TM\ox TM \, .
\edm
On the other side, a connection $\nabla$ is metric if and only if
and only if $A$ belongs to the space
\bdm
\mathcal{A}^g\ :=\ TM\ox\wedge^2 TM \ = \ \{A \in\ox^3 TM \ | 
\ A(X,V,W) +  A(X,W,V) \ =\ 0\} \, .
\edm
The real orthogonal group  $\Orth(n,\R)$ acts on both tensor representations
$\mathcal{T}$ and $\mathcal{A}^g$
in a natural way by $g\cdot\T(X,Y,Z):=\T(g^{-1}X,g^{-1}Y,g^{-1}Z)$
for $g\in\Orth(n,\R)$. 
\begin{prop}\label{classes}
For $n\geq 3$, the space $\mathcal{T}$ of possible torsion tensors
splits under $\Orth(n,\R)$ into the sum of three irreducible representations,
$\mathcal{T}\cong TM\op \wedge^3 T M \op \mathcal{T}'$, as
does $\mathcal{A}^g$.
Furthermore, an equivariant bijection 
$\Phi:\ \mathcal{A}^g \ra \mathcal{T}$ is given by
($A\in \mathcal{A}^g, \T \in\mathcal{T} $)
\begin{eqnarray*} 
\Phi(A)(X,Y,Z) & = & A(X,Y,Z)-A(Y,X,Z)\, ,\\ 
2\Phi^{-1}(\T)(X,Y,Z) & = & \T(X,Y,Z)-\T(Y,Z,X)+\T(Z,X,Y) \ .
\end{eqnarray*} 
The map $\Phi$ is a multiple of the identity precisely on $\wedge^3 T M$.
\end{prop}
\begin{proof}
It is clear that $\mathcal{T}$ and $\mathcal{A}^g$ split into the
same irreducible summands under $\Orth(n,\R)$. Hence, we concentrate
on $\mathcal{T}$.
There exist  two $\Orth(n,\R)$-equivariant 
contractions from $\mathcal{T}$ into irreducible 
$\Orth(n,\R)$-representations,
\bdm  \Phi_1:\ \mathcal{T}\lra \wedge^3 TM,\quad
\Phi_2: \ \mathcal{T} \lra TM,
\edm
given by 
\bdm
\Phi_1(\T)\ =\ \frac{1}{3}\,\mathfrak{S}\, \T(X,Y,Z),\quad
\Phi_2(\T)\ =\ \sum_{i=1}^n \T(e_{i+1},e_i,e_{i+1})e_i.
 \edm
Here, $\mathfrak{S}$ denotes antisymmetrisation with respect to all 
arguments and $e_1,\ldots,e_n$ is any orthonormal basis of $T M$. 
Vice versa, $TM$ can be realized as an irreducible subspace of 
$\mathcal{T}$ via $\Phi_2^{-1}:T M\ra \mathcal{T}$,
\bdm
V\mapsto \T_V,\, \T_V(X,Y,Z):=g( X,Z) g( V, Y) -g(Y,Z)g(V, X) .
\edm
All in all, we identified two irreducible summands of $\mathcal{T}$,
$\wedge^3 TM\subset\ker \Phi_2$ and $TM\subset\ker\Phi_1$. 
A dimensional argument
shows that $\mathcal{T}':=\ker\Phi_1\cap\ker\Phi_2$ is not empty.
In fact, one easily checks that it is irreducible under the action
of $\Orth(n,\R)$, and a routine calculation proves all claims about the
isomorphism $\Phi$.
\end{proof}
\noindent
The eight classes of linear connections are now defined by the possible
parts of their torsions $\T$ in these components.  If one looks at the class
of linear \emph{metric} connections, then these are also uniquely determined
by their torsion, since $\Phi^{-1}$ reconstructs $A$ from $\T$.
For general connections, $\T$ determines $A$ only up to a 
contribution from the complement of $\mathcal{A}^g$ inside
$\ox^3 TM$, that is, from $TM\ox S^2 TM$. Since this space
splits itself into two irreducible subspaces, one might as well speak
of a total of $16$ classes in the general situation. The nice lecture notes 
by Tricerri and Vanhecke \cite{TriVan} use a similar approach in order
to classify homogeneous spaces by the algebraic properties of the torsion
of the canonical connection. They construct homogeneous examples of all
classes, and study their ``richness''.

\smallskip\noindent
The described decompositions shows that a natural class of
non-torsion free metric connections are those with skew-symmetric
torsion form. We obtain a geometric characterization of 
these connections.  
\begin{cor}\label{charact-skewI}
A connection $\nabla$ on M is metric and geodesics preserving precisely
if its torsion $\T$ lies in $\wedge^3 TM$.
In this case, $2\cdot A=\T$ holds,
\bdm
\nabla_X Y \ = \ \nabla^g_X Y \, + \, \frac{1}{2} \cdot \T(X,Y,-),
\edm 
and the $\nabla$-Killing vector fields coincide with the Riemannian Killing 
vector fields. 
\end{cor}
\begin{proof}
If $\nabla$ preserves geodesics, $2\cdot A=\T$ by Lemma~\ref{same-geodesics}.
If $\nabla$ is also metric, $A$ needs in addition to lie in
the component of $\mathcal{A}^g$ that yields a torsion proportional to
$A$, which is $\wedge^3 TM$ by Proposition~\ref{classes}.
\end{proof}
%
\section{The holonomy of spinor connections with constant torsion in $\R^n$}
\noindent
%
We consider the euclidian vector space $\R^n$ equipped with its standard
inner product. The exterior algebra $\Lambda^*(\R^n)$ and the Clifford algebra
$\mathrm{Cl}(\R^n)$ are -- treated as vector spaces only -- equivalent 
$\SO(n)$-representations. Denote by $\Delta_n$ the complex vector space 
of all $n$-dimensional spinors.
The Clifford algebra - and henceforth the exterior algebra, too - 
acts on $\Delta_n$. We denote by $\T \cdot 
\psi$ the corresponding action of a $k$-form $\T$ on a spinor
$\psi \in \Delta_n$. It is $\SO(n)$-equivariant and called the {\it Clifford 
multiplication} of a spinor by a $k$-form. The Clifford 
algebra is an associative algebra and there is an underlying 
Lie algebra structure,
\bdm
[ \alpha, \beta ] \ = \ \alpha \cdot \beta \, - \, \beta \cdot \alpha , 
\quad \alpha, \beta \in \mathrm{Cl}(\R^n) \, .
\edm
We denote the corresponding Lie algebra by $\mathfrak{cl}(\R^n)$. The 
Lie algebra $\mathfrak{so}(n)$ of the special orthogonal group is a subalgebra
of $\mathfrak{cl}(\R^n)$, 
\bdm
\mathfrak{so}(n) \ = \ \mathrm{Lin} \big\{X \cdot Y  : X,
Y \in \R^n \ \mathrm{and} 
\  \langle X  , \, Y \rangle \ = \ 0 \big\} \, \subset \, 
\mathfrak{cl}(\R^n) \, .
\edm
Consider an algebraic $k$-form $\T \in \Lambda^k(\R^n)$ and 
denote by $\mathrm{G}_{\T}$ the group of all orthogonal
transformation of $\R^n$ preserving the form $\T$. Let 
$\mathfrak{g}_{\T}$ be its Lie algebra. We associate with any 
exterior form its covariant derivative $\nabla^{\T}$ acting on spinor 
fields $\psi : \R^n \rightarrow \Delta_n$ by the formula
\bdm
\nabla^{\T}_X \psi \ := \ \nabla^{g}_X \psi \, + \,  
(X \haken \T) \cdot \psi \, .
\edm 
Here, $\nabla^{g}$ denotes the Levi-Civita connection. 
For a $3$-form $\T \in \Lambda^3(\R^n)$, the spinorial covariant 
derivative $\nabla^{\T}$ is induced by a linear metric 
connection with torsion tensor $2 \cdot \T$,
\bdm
\nabla^{\T}_X Y \ := \ \nabla^{g}_X Y \, + \, 2 
\cdot \T(X,Y,-) \, .
\edm
 For a general exterior form $\T$, we introduce
a new Lie algebra $\mathfrak{g}_{\T}^*$ that is a subalgebra of
$\mathfrak{cl}(\R^n)$.
\begin{dfn}
Let $\T$ be an exterior form on $\R^n$. The Lie algebra 
$\mathfrak{g}_{\T}^*$ is the subalgebra of $\mathfrak{cl}(\R^n)$ generated by all elements $X \haken \T$, where $X \in \R^n$ is a vector.
\end{dfn}
\noindent
The Lie algebra $\mathfrak{g}_{\T}^*$ is invariant under
the action of the isotropy group $\mathrm{G}_{\T}$. The derived
algebra $\big[\mathfrak{g}_{\T}^* \, , \, \mathfrak{g}_{\T}^*
\big]$ is the Lie algebra generated by all curvature transformations
of the spinorial connection $\nabla^{\T}$. It is the Lie algebra 
of the infinitesimal holonomy group of the spinorial covariant 
derivative $\nabla^{\T}$ (see \cite{KN}, Chapter II, Section 10):
\begin{dfn} 
Let $\T$ be an exterior form on $\R^n$. The Lie algebra
\bdm
\mathfrak{h}_{\T}^* \ := \ \big[\mathfrak{g}_{\T}^* , \, 
\mathfrak{g}_{\T}^*\big] \ \subset \ \mathfrak{cl}(\R^n)  
\edm
is called the \emph{infinitesimal holonomy algebra} of the exterior form $\T$.
\end{dfn} 
\noindent
The Lie algebra $\mathfrak{h}_{\T}^*$ is invariant under
the action of the isotropy group $\mathrm{G}_{\T}$, too.
%
For a $3$-form $\T$, the Lie algebras 
$\mathfrak{g}_{\T}^* , \, \mathfrak{h}_{\T}^* \subset \mathfrak{so}(n)$ 
are subalgebras of the Lie algebra of the orthogonal group. 
This inclusion reflects again the fact that the corresponding 
spinor derivative $\nabla^{\T}$ is induced by a linear metric
connection. The following proposition generalizes this observation.
\begin{prop} 
If $\T $ is a $k$-form with
$k + \binom{k-1}{2} \equiv 0 \ \mathrm{mod}\, 2$,
then $\g_{\T}^*$ is a compact Lie algebra.
\end{prop}
\begin{proof} We consider the complex spin representation of the 
Clifford algebra. There exists a hermitian product on $\Delta_n$ such
that
\bdm
\big( X \cdot \psi \, , \, \psi_1 \big) \ + \ \big( \psi \, , \, X \cdot 
\psi_1 \big) \ = \ 0
\edm
for all vectors $X \in \R^n$ and all spinors $\psi,\psi_1 \in \Delta_n$.
Then, under the condition for the degree of the form $\T$, all
endomorphisms $X \haken \T$ acting on $\Delta_n$ are skew-symmetric.
\end{proof}
\noindent
The following proposition is a special case of the general holonomy theory.
For completeness, let us sketch its proof.
\begin{prop} 
There exists a non-trivial $\nabla^{\T}$-parallel spinor field $\psi : \R^n 
\rightarrow \Delta_n$, 
\bdm
\nabla^{\T}_X \psi \ = \ X(\psi) \, + \, (X \haken \T) \cdot
\psi \ = \ 0,
\edm
if and only if there exists a constant spinor $\psi_0 \in \Delta_n$ such that
$\mathfrak{h}_{\T}^* \cdot \psi_0 =  0$.
\end{prop} 
\begin{proof} 
If $\psi : \R^n \rightarrow \Delta_n$ is 
$\nabla^{\T}$-parallel, we differentiate it twice with respect to
arbitrary vectors $X,Y \in \R^n$. Then we  obtain the condition
\bdm
\big[ X \haken \T \, , \, Y \haken \T \big] \cdot \psi \ 
= \ 0 \, ,
\edm
i.e., $\mathfrak{h}_{\T}^* \cdot \psi = 0$. Conversely, if $\psi_0 
\in \Delta_0$ is a spinor such that $\mathfrak{h}_{\T}^* \cdot \psi_0
= 0$, we define the spinor field $\psi : \R^n \rightarrow \Delta_n$ by the 
formula
\bdm
\psi(m) \ := \ \mathrm{Exp}( - m \haken \T) \cdot \psi_0 \, , \quad 
m \in \R^n \ .
\edm
An easy computation yields that $X(\psi)(m) +  (X \haken \T) \cdot 
\psi(m)$ is given by the formula
\bdm
\mathrm{Ad}\Big(\mathrm{Exp}(m \haken \T)\Big)\Big( \frac{[m \haken \T , X \haken \T]}{2} \, + \, 
 \frac{[ m \haken \T \, , \, [m \haken \T , X \haken 
\T]]}{6} \, + \cdots \Big) \cdot \psi_0 \ .
\edm 
The commutators $[m \haken \T, X \haken \T] \, $ etc. are 
in $\mathfrak{h}_{\T}^*$ and the adjoint action $\mathrm{Ad}(\mathrm{Exp}(m \haken \T))$ preserves the holonomy algebra $\mathfrak{h}_{\T}^*$ since 
$m \haken \T \in \mathfrak{g}_{\T}^*$. 
\end{proof}
\begin{cor}\label{parallelspinor}  
Let $\T$ be an exterior form such that the Lie algebra
$\mathfrak{g}_{\T}^*$ is perfect, $\mathfrak{h}_{\T}^* 
= \mathfrak{g}_{\T}^*$. Then any $\nabla^{\T}$-parallel 
spinor field $\psi : \R^n \rightarrow \Delta_n$ , 
\bdm
\nabla^{\T}_X \psi \ = \ X(\psi) \, + \, (X \haken \T) \cdot
\psi \ = \ 0 \, ,
\edm
is constant and $\mathfrak{g}_{\T}^* \cdot \psi = 0$.
\end{cor}
\begin{proof} 
Any parallel spinor field satisfies the condition 
$\mathfrak{h}_{\T}^* \cdot \psi = 0$. By assumption, we obtain
$\mathfrak{g}_{\T}^* \cdot \psi = 0$ and the differential equation 
yields $X(\psi) = 0$, i.e., $\psi$
is constant. 
\end{proof}
\begin{exa} 
If $\T \in \Lambda^1(\R^n)$ is a $1$-form, 
the Lie algebra $\mathfrak{g}_{\T}^*$ is generated by one
element $1 \in \mathfrak{cl}(\R^n)$ and $\mathfrak{g}_{\T}^* 
=  \R$, $\mathfrak{h}_{\T}^* =0$. The general solution of the
equation $\nabla^{\T} \psi = 0$ is
\bdm
\psi(m) \ = \ e^{ -  \langle m , \T \rangle } \cdot \psi_0 \, , 
\quad m \in \R^n \, ,
\edm
where $\psi_0$ is constant.
\end{exa}
\noindent
We denote by $e_1, \ldots , e_n$ an orthonormal frame on $\R^n$, and 
abbreviate as
$e_{ijk\ldots}$ the exterior product $e_i \wedge e_j \wedge e_k \wedge 
\ldots$ of  $1$-forms.
\begin{exa} 
Any $2$-form $\T \in \Lambda^2(\R^n)$ of rank $2k$ is
equivalent to $A_1 \cdot e_{12} + \cdots + A_k \cdot e_{2k-1,2k}$. The 
Lie algebra $\mathfrak{g}_{\T}^*$ is generated by
the elements $e_1, e_2, \cdots , e_{2k-1}, e_{2k}$. It is isomorphic
to the Lie algebra $\spin(2k+1)$. In particular, if $n = 8$ then 
$\Delta_8 = \R^{16}$ is a real, $16$-dimensional and the spinorial holonomy
algebra of a generic $2$-form in eight variables is the unique $16$-dimensional
irreducible representation of $\spin(9)$.

\end{exa}
\begin{exa}
Consider the $4$-form $\mathrm{T} = e_{1234} + e_{3456} \in \Lambda^4(\R^6)$.
The Clifford algebra $\mathrm{Cl}(\R^6) = \mathrm{End}(\R^8)$ is isomorphic
to the algebra of all endomorphisms of an $8$-dimensional real vector
space and $\mathfrak{g}_{\mathrm{T}}^*$ is the Lie algebra generated by the
elements
\bdm
e_{234} , \quad e_{134} , \quad e_{124}+e_{456}, \quad e_{123} + e_{356}, 
\quad e_{346} , \quad e_{345} .
\edm 
A computation of the whole Lie algebra yields the result that 
$\mathfrak{g}_{\mathrm{T}}^*$ is isomorphic  to the Lie algebra
$\mathfrak{e}(6)$ of the euclidian group. 
\end{exa} 
\begin{exa} 
Consider the volume form $\mathrm{T} = e_{123456}$ in $\R^6$.
The subalgebra $\g_{\T}^*$ of $\mathrm{Cl}(\R^6) = \mathrm{End}(\R^8)$
is isomorphic to the compact Lie algebra $\spin(7)$. Indeed, it is 
generated by the Lie algebra $\spin(6)$ and
all elements of degree five.
\end{exa}
\begin{exa} Let us discuss the holonomy algebra of a more complicated
$4$-form in seven variables,
\bdm
\T \ = \ e_{12}\cdot (e_{34} \, - \, e_{56}) \, - \,
e_{17}\cdot (e_{45} \, - \, e_{36}) \, - \,
e_{27}\cdot (e_{35} \, + \, e_{46}) \, - \,
e_{3456} \ . 
\edm
The $7$-dimensional spin representation is real and we describe 
the holonomy algebra $\g_{\T}^*$ using the spin representation
$\mathfrak{cl}(\R^7) \rightarrow \mathfrak{gl}(\Delta_7) = 
\mathfrak{gl}(\R^8)$ of the Clifford algebra. For this purpose, we introduce
the matrices
\bdm
A_1:=\begin{bmatrix} 0 & 0 & 0 & 0 \\ 0 & 0 & 0 & 0 \\ 1 & 0 & 0 & 0\\
0 & 0 & 0 & 0\end{bmatrix},\, 
A_2:=\begin{bmatrix} 0 & 0 & 0 & 0 \\ 0 & 0 & 0 & 0 \\ 0 & 1 & 0 & 0\\
0 & 0 & 0 & 0\end{bmatrix},\,  
A_3:=\begin{bmatrix} 0 & 0 & 0 & 0 \\ 0 & 0 & 0 & 0 \\ 0 & 0 & 1 & 0\\
0 & 0 & 0 & 0\end{bmatrix},\, 
A_4:=\begin{bmatrix} 0 & 0 & 0 & 0 \\ 0 & 0 & 0 & 0 \\ 0 & 0 & 0 & 1\\
0 & 0 & 0 & 0\end{bmatrix} . 
\edm
The holonomy algebra, treated as a subalgebra of $\mathfrak{gl}(\R^8)$,
is the Lie algebra generated by the following seven matrices:
\bdm
B_1:=\begin{bmatrix}0 & A_1 \\ A_1^t & 0\end{bmatrix} ,\ 
B_2:=\begin{bmatrix}0 & A_2 \\ A_2^t & 0\end{bmatrix} ,\
B_3:=\begin{bmatrix}0 & A_3 \\ A_3^t & 0\end{bmatrix} , \
B_4:=\begin{bmatrix}0 & A_4 \\ A_4^t & 0\end{bmatrix} ,
\edm

\bdm
B_5:=\begin{bmatrix}A_1 + A_1^t & 0 \\ 0 & 0\end{bmatrix} ,\ \
B_6:=\begin{bmatrix}A_2 + A_2^t & 0 \\ 0 & 0\end{bmatrix} ,\ \
B_7:=\begin{bmatrix}A_4 + A_4^t & 0 \\ 0 & 0\end{bmatrix} .  
\edm
An investigation of the commutators of these matrices yields the result
that $\mathfrak{g}_{\T}^*$ is a $46$-dimensional subalgebra
of $\mathfrak{gl}(\R^8)$,
\bdm
\mathfrak{g}_{\T}^* \ = \ \Big\{ \begin{bmatrix}X & A \\ A^t & Y\end{bmatrix} 
\,  : \ X , Y \in \mathfrak{sl}(\R^4) \ \mathrm{and} \ A \in 
\mathfrak{gl}(\R^4) \Big\} \ .
\edm
No spinor is fixed by the holonomy group of the connection $\nabla^{\T}$,
i.e., in the flat space $\nabla^{\T}$-parallel spinors do not exist. Later
we will see that this torsion form occurs in certain compact Riemannian
manifolds in a natural way. On these non flat spaces there exist 
$\nabla^{\T}$-parallel spinors, see Theorem \ref{3Sas4form}.
\end{exa}
\section{Constant $3$-forms in $\R^n$ and their holonomy algebra}\noindent
%
We will study $3$-forms $\T \in \Lambda^3(\R^n)$ and their
Lie algebras $\mathfrak{g}_{ \T}^*$. To begin with, let us consider
some examples.
\begin{exa}
This is the place to discuss Cartan's first example of a space
with torsion (see \cite[p.~595]{Cartan1}). Consider $\R^3$ with its usual
euclidian metric, and the connection
\bdm
\nabla_X Y\ =\ \nabla^g_X Y - X\x Y,
\edm
corresponding, of course, to the choice $\T=-2 \cdot e_1\wedge e_2\wedge e_3$.
Cartan observed correctly that this connection has same geodesics than
$\nabla^g$, but induces a different parallel transport\footnote{
``Deux tri\`edres [\ldots] de $\mathcal{E}$ seront parall\`eles lorsque les 
tri\`edres
correspondants de E [l'espace euclidien classique] pourront se d\'eduire
l'un de l'autre par un d\'eplacement h\'elico\"{\i}dal de pas donn\'e, de sens
donn\'e[\ldots]. L'espace $\mathcal{E}$ ainsi d\'efini admet un groupe de 
transformations \`a 6 param\`etres : ce serait notre espace ordinaire vu par
des observateurs dont toutes les perceptions seraient tordues.'' loc.cit.}.
Indeed, consider the $z$-axis $\gamma(t)=(0,0,t)$, a geodesic, and the
vector field $V$ which, in every point $\gamma(t)$, consists of the
vector $(\cos t,\sin t,0)$. Then one checks immediately that
$\nabla^g_{\dot{\gamma}}V= \dot{\gamma}\x V$, that is, the vector
$V$ is parallel transported according to a helicoidal movement.
If we now transport the vector along the edges of a closed triangle,
it will be rotated around three linearly independent axes, hence
the holonomy algebra is $\mathfrak{g}^*_{\T}=\mathfrak{h}^*_{\T}=
\mathfrak{so}(3)$.
\end{exa}
\begin{exa} 
Any $3$-form in $\R^4$ is equivalent to one of the forms
$\T =  a \cdot e_{123}$, hence the same argument as in the previous
example yields that $\mathfrak{g}_{\T}^* = 0$ or $\mathfrak{so}(3)$.
\end{exa}
\begin{exa} 
Any $3$-form in $\R^5$ is equivalent to one of the forms
$\T =  a \cdot e_{123} \, + \, b \cdot e_{345}$. The corresponding algebras
are $\mathfrak{g}_{\T} =  \mathfrak{so}(5), \ \mathfrak{so}(3)\oplus 
\mathfrak{so}(2), 
\ 0 \ $ and  $\mathfrak{g}_{\T}^* = 0 , \ \mathfrak{so}(3), 
\ \mathfrak{so}(5)$.
\end{exa}
\begin{exa} 
In $\R^7$, we consider the $3$-form
$\T =  e_{127} + e_{135} - e_{146} - e_{236} - e_{245} + e_{347} + 
e_{567}$.  
Its isotropy algebra $\mathfrak{g}_{\T}$ is isomorphic to the
exceptional Lie algebra $\mathfrak{g}_2$. Moreover, $\mathfrak{so}(7)$
splits into two $\mathrm{G}_2$-irreducible components,
$\mathfrak{so}(7) =  \mathfrak{g}_{\T} \oplus \mathfrak{m}$. The orthogonal
complement $\mathfrak{m}$ of $\mathfrak{g}_{\T}$
 coincides with the space of all inner products $X \haken \T$. 
The Lie algebra generated by these elements 
is isomorphic to $\mathfrak{so}(7)$.
To summarize, we obtain
$\mathfrak{g}_{\T} =  \mathfrak{g}_2$  and 
$\mathfrak{g}_{\T}^* =  \mathfrak{so}(7)$. 
\end{exa}
\noindent
The first Proposition estimates the dimension of the Lie algebra
$\mathfrak{g}_{\T}^*$ from below.
\begin{prop}\label{inequ}  
Let $\T \in \Lambda^3(\R^n)$ be a $3$-form 
and $\Phi_{\T} : \R^n \rightarrow \mathfrak{g}_{\T}^*$
be the map defined by the formula $\Phi_{\T}(X) := X \haken 
\T$. Then $\T$ depends only on the orthogonal
complement $\mathrm{Ker}(\Phi_{\T})^{\perp}$, 
\bdm
\T \ \in \ \Lambda^3(\mathrm{Ker}(\Phi_{\T})^{\perp}) \ .
\edm
In particular, if $\T$ is a $3$-form which can not be reduced to
a lower dimensional subspace, then
\bdm
n \ \leq \ \mathrm{dim}(\mathfrak{g}_{\T}^*) \, .
\edm
\end{prop} 
\noindent
Next, we investigate the representation of the Lie algebra 
$\mathfrak{g}_{\T}^*$ in $\R^n$.
\begin{prop} \label{3formsplitting} 
The representation $(\mathfrak{g}_{\T}^* \, , \, \R^n)$ is reducible if and only if there exist a proper subspace $\mathrm{V} \subset 
\R^n$ and two $3$-forms $\T_1 \in \Lambda^3(\mathrm{V})$ and 
 $\T_2 \in \Lambda^3(\mathrm{V}^{\perp})$ such that $\T = 
\T_1 + \T_2$. In this case, the Lie algebra 
$\mathfrak{g}_{\T}^*$ decomposes into
\bdm
\mathfrak{g}_{\T}^*  \ = \ \mathfrak{g}_{\mathrm{T_1}}^* \oplus 
\mathfrak{g}_{\mathrm{T_2}}^*  . 
\edm
\end{prop}
\begin{proof} 
Consider a $\mathfrak{g}_{\T}^*$-invariant subspace
$\mathrm{V} \subset \R^n$ and fix a basis $e_1, \cdots , e_k$ in $\mathrm{V}$
as well as a basis $e_{k+1}, \cdots , e_n$ in its orthogonal complement
$\mathrm{V}^{\perp}$. Then, for any vector $X \in \R^n$, and any pair of
indices $1 \leq i \leq k$, $k+1 \leq \alpha \leq n$, we obtain
\bdm
\T(X, \, e_i, \, e_{\alpha}) \ = \ 0 \, .
\edm
Since $\T$ is skew-symmetric, we conclude
\bdm
\T(e_i, \, e_j, \, e_{\alpha}) \ = \ 0 , \quad \mathrm{and} \quad
\T(e_i, \, e_{\alpha}, \, e_{\beta}) \ = \ 0 \, . \qedhere
\edm
\end{proof}
\noindent
The following Proposition restricts the type of the Lie algebra 
$\mathfrak{g}_{\T}^*$. In particular, it cannot be contained in the 
Lie algebra 
$\mathfrak{u}(k) \subset \mathfrak{so}(2k)$ of the unitary group.
\begin{prop}\label{invform}  
Let $\T$ be a $3$-form in $\R^{2k}$ and suppose that there
exists a $2$-form $\Omega$ such that
\bdm
\Omega^{k} \ \neq \ 0 \quad \mathrm{and} \quad [\, \mathfrak{g}_{\T}^*, \, \Omega \, ] \ = \ 0 \, .
\edm
Then $\T$ is zero, $\T = 0$.
\end{prop} 
\begin{proof} 
We fix an orthonormal basis in $\R^{2k}$ such that the $2$-form 
$\Omega$ is given by
\bdm
\Omega \ = \ A_1 \cdot e_{12} \, + \cdots \, + \, A_k \cdot e_{2k-1,2k} \, , 
\quad A_1 \cdot \ldots \cdot A_k \ \neq \ 0 \, . 
\edm
The condition 
$[\mathfrak{g}_{\T}^* , \, \Omega] = 0$ is equivalent to the equations
\bdm
\sum_{j=1}^{2k} \Omega_{\alpha j} \cdot \T_{\beta j \gamma} \ = \ 
\sum_{j=1}^{2k} \T_{\beta \alpha j} \cdot \Omega_{j \gamma}
\edm
for any triple $1 \leq \alpha, \beta, \gamma \leq 2k$. Using the special form
of $\Omega$ we obtain the equations ($1 \leq \alpha, \gamma \leq k$):
\bdm
A_{\alpha} \cdot \T_{\beta, 2 \alpha, 2 \gamma -1} \ = \ - \, 
A_{\gamma} \cdot \T_{\beta, 2 \alpha -1, 2 \gamma}
\edm
and 
\bdm
A_{\alpha} \cdot \T_{\beta, 2 \alpha-1, 2 \gamma -1} \ = \  
A_{\gamma} \cdot \T_{\beta, 2 \alpha, 2 \gamma} \, .
\edm
The latter system of algebraic equations implies that $\T = 0$ 
vanishes. Indeed, let us compute -- for example -- 
$\T_{\beta, 2 \alpha, 2 \gamma -1}$. In case $\beta$ is odd, we have
\begin{eqnarray*}
A_{\alpha} \cdot \T_{\beta, 2 \alpha, 2 \gamma -1} & = & 
- \, A_{\gamma} \cdot \T_{\beta, 2 \alpha - 1, 
2 \gamma} \ = \ 
A_{\gamma} \cdot \T_{2 \alpha - 1, \beta, 2 \gamma} \ = \ - \, 
A_{(\beta + 1)/2} \cdot \T_{2 \gamma -1, 
2 \alpha - 1, \beta + 1}\\
 & = & A_{\alpha} \cdot \T_{2 \gamma - 1, 2 \alpha, \beta} \ = 
\ - \, A_{\alpha} \cdot \T_{\beta, 2 \alpha, 2 \gamma -1} \ .
\end{eqnarray*}
In case $\beta$ is even, a similar computations yields the formula
\bdm
\big[A_{\alpha}\big]^2 \cdot \T_{\beta, 2 \alpha, 2 \gamma -1} \ = \ - \big[A_{\beta/2}\big]^2 \cdot \T_{\beta, 2 \alpha, 2 \gamma -1} \ . \qedhere
\edm
\end{proof}
\begin{thm} \label{3formhalbeinfach} 
For any $3$-form $\T \in \Lambda^3(\R^n)$, the Lie algebra
$\mathfrak{g}_{\T}^*$ is semisimple and coincides with the holonomy
algebra $\mathfrak{h}_{\T}^*$.
\end{thm}  
\begin{proof} 
According to Proposition \ref{3formsplitting} we assume that the
representation $(\mathfrak{g}_{\T}^* \, , \, \R^n)$ is irreducible.
The Lie algebra $\mathfrak{g}_{\T}^*$ splits into the 
holonomy algebra $\mathfrak{h}_{\T}^*$ and the
center $\mathfrak{z}(\mathfrak{g}_{\T}^*)$. Suppose that the
center $\mathfrak{z}$ is non trivial, i.e., that there exist a 
$2$-form $\Omega$ such that
\bdm
\big[\mathfrak{g}_{\T}^* \, , \, \Omega \big] \ = \ 0 \, . 
\edm  
We split the euclidian space into
\bdm
\R^n \ = \ \mathrm{Ker}(\Omega) \, \oplus \, \mathrm{Ker}(\Omega)^{\perp} 
\edm
and observe that both subspaces are $\mathfrak{g}_{\T}^*$-invariant.
Since $\mathrm{Ker}(\Omega) \neq 0$ and the representation 
$(\mathfrak{g}_{\T}^* \, , \, \R^n)$ is irreducible, we conclude that
$\mathrm{Ker}(\Omega) = 0$. In particular, the dimension $n = 2k$ is even and
$\Omega^k \neq 0$. Finally, we obtain $\T = 0$ by Proposition~\ref{invform}.  
\end{proof}
\noindent
A second restriction for the algebra $\mathfrak{g}_{\T}^*$ results
from the observation that it is not contained in the isotropy Lie algebra
of a spinor. This fact implies that there are no 
$\nabla^{\T}$-parallel spinors in $\R^n$ for $\T \neq 0 $.
Furthermore, certain semisimple Lie groups 
cannot occur as holonomy groups of $3$-form in $\R^n$. 
In dimensions $n \leq 9$, 
where the group $\Spin(n)$ acts transitively
on the set of spinors of length one, the proof is a consequence
of a direct
algebraic computation. 
For example, in dimension $n = 8$, a general $3$-form depends on $56$ 
parameters
and $\mathfrak{g}_{\T}^* \cdot \psi = 0$ is a system consisting again
of at least $56$ linear equations. In higher dimensions, we have
to avoid the problem of the unknown orbit structure of the spin representation.
We use a global argument here, but it would be interesting to find 
a purely algebraic proof. 
\begin{thm}\label{pareucl} 
Let $\T \in \Lambda^3(\R^n)$ be a $3$-form. If there exists a non trivial 
spinor $\psi \in \Delta_n$ such that $\mathfrak{g}_{\T}^* \cdot \psi = 0$, 
then $\T = 0$.
\end{thm} 
\begin{proof} Consider the compact, flat torus $\R^n / \Z^n$. Since 
$\T$ and $\psi \in \Delta_n$ are constant, both are geometric objects on
the torus. In particular, with respect to the trivial spin structure
of the torus, $\psi$ is a $\nabla^{\T}$-parallel spinor field on
$\R^n / \Z^n$. The integral formula of Theorem \ref{parallelflach2} yields that
$\T = 0$. 
\end{proof} 
\begin{cor}\label{parallelflach} 
Let $\T \in \Lambda^3(\R^n)$ a $3$-form. If there exists a non trivial 
solution $\psi : \R^n \rightarrow \Delta_n$ of the equation
\bdm
\nabla^{\T}_X \psi \ = \ X(\psi) \, + \, (X \haken \T) \cdot
\psi \ = \ 0 ,
\edm
then $\T = 0$ and $\psi$ is constant.
\end{cor}
\begin{proof} 
Suppose that a non trivial parallel spinor exists. By
Corollary \ref{parallelspinor} and Theorem \ref{3formhalbeinfach}, 
we conclude that $\psi$ is constant 
and $\mathfrak{g}_{\T}^* \cdot \psi = 0$. Theorem \ref{pareucl} yields now
that the $3$-form $\T = 0$ vanishes.
\end{proof}  
\noindent
In low dimensions, we obtain a complete list
of all possible holonomy algebras $\mathfrak{h}_{\T}^*$:
\begin{itemize}
\item $n \ =\ 5 :  \quad \mathfrak{h}_{\T}^* \ = \ 0  , \ 
\mathfrak{so}(3) , \ \mathfrak{so}(5) $. 
\item $n \ =\ 6  :\quad \mathfrak{h}_{\T}^* \ = \ 0  , \ 
\mathfrak{so}(3) , \ \mathfrak{so}(5) , \ \mathfrak{so}(3) \, \oplus \, 
\mathfrak{so}(3) , \ \mathfrak{so}(6) $.  
\item $n \ =\ 7  : \quad \mathfrak{h}_{\T}^* \ = \ 0 , \ 
\mathfrak{so}(3) , \ \mathfrak{so}(5) , \ \mathfrak{so}(3) \, \oplus \, 
\mathfrak{so}(3) , \ \mathfrak{so}(6) , \ \mathfrak{so}(7) $. 
\end{itemize}
\noindent
Starting from dimension eight, there occur representations
of all semisimple Lie algebras as the holonomy algebra of certain
$3$-form in euclidian space. Indeed,  
suppose that
the euclidian space $\R^n = \g$ is a compact Lie algebra, and the inner
product and the Lie bracket are related by the condition
\bdm
\langle\, [X \, , \, Y ] \, , \, Z \, \rangle \, + \, \langle \, 
Y \, , \, [X \, , \, Z ] \, \rangle \ = \ 0 \, .
\edm
Then $\T(X,Y,Z):= \langle [X,Y],Z \rangle$ is a $3$-form in $\R^n = \g$ and
we obtain
\bdm
X \haken \T \ = \ \mathrm{ad}(X) \ \in \ \so(\g) \ = 
\ \so(n) \, .
\edm
The Lie algebra $\g^*_{\T}$ is the image of the Lie algebra
$\g$ under the adjoint representation.
Consequently, we have a series of representations occurring
for some $3$-form. 
\begin{prop} 
The adjoint representation of any compact, semisimple
Lie algebra $\g$ is the holonomy algebra of a certain $3$-form
with constant coefficients in euclidian space $\g = \R^n$.
\end{prop}
\noindent   
The first
interesting example is the $8$-dimensional Lie algebra $\su(3)$. It
yields a $3$-form in $\R^8$ such that $\g^*_{\T} = \su(3)$ and the
inclusion $\su(3) \subset \so(8)$ is the adjoint representation. This example
realizes the lower bound in the dimension estimate of Proposition \ref{inequ}.
The latter series of examples generalizes to Riemannian naturally 
reductive spaces
$G/H$. Decompose the Lie algebra
\bdm
\g \ = \ \h \, \oplus \, \m  , \quad \mathrm{Ad}(H)(\m) \, \subset \, \m \, ,
\edm 
and consider the canonical connection of the reductive space. Its torsion
form is given by the formula
\bdm
\T(X,Y,Z) \ = \ - \langle \, [X, Y]_{\m}\, , \, Z \rangle \, , 
\quad X,Y,Z \in \m \, ,
\edm
where $[ \ , \ ]_{\m}$ denotes the $\m$-part of the Lie bracket. Consider
the euclidian space $\m$ and the $3$-form $\T$. Then $\g_{\T}^*$
is the Lie subalgebra of $\so(\m)$ generated  
by the subspace $\m \rightarrow \so(\m)$, where this map is given by the
formula
\bdm
Z \ \longrightarrow \ Z \haken \T \, , 
\quad (Z \haken \T)(X) \ = \ [X  ,  Z]_{\m} , \quad Z \in \m \, .
\edm
In general, this is {\it not} the isotropy representation of the reductive 
space, but related to the holonomy of its Levi-Civita connection (see
\cite{KN2}). 

\smallskip\noindent
Let us discuss the question which irreducible representations
$(\g, \,\R^n)$ of a semisimple Lie algebra $\g$ can occur
for a $3$-form. We already know some restrictions. In even dimensions,
the $\g$-action cannot preserve a non-degenerate $2$-form and, in any
dimension, the lift into the spin representation cannot preserve
a spinor. In order to formulate a further restriction we introduce --
in analogy to the prolongation of a linear Lie algebra 
(see \cite[note 13]{KN2}) --
an \emph{antisymmetric prolongation} of a representation of a compact
semisimple Lie algebra by
\bdm
\T(\g , \R^n) \ := \ \big\{ \T \in \Lambda^3(\R^n) \, : \, X \haken \T 
\in \g \ \ \text{for  any} \ X \in \R^n \ \big\} .
\edm 
The subspace $\T(\g,\R^n) \subset \Lambda^3(\R^n)$ is $\g$-invariant.
A $3$-form $\T$  belongs to this space if and only if its Lie algebra is
contained in $\g_{\T}^* \subset \g$. In particular, we can 
formulate a necessary condition.
\begin{prop} If a representation $(\g,\R^n)$ of a compact, semisimple
Lie algebra is realized by some $3$-form $\T \in \Lambda^3(\R^n)$, then
$\T(\g , \R^n) \not= 0$ is non trivial.
\end{prop}
\begin{exa} 
The unique irreducible $16$-dimensional 
representation $\spin(9) \subset \so(16)$ of the Lie algebra $\spin(9)$
does not admit invariant, non degenerate $2$-forms in $\R^{16}$ or 
invariant spinors in $\Delta_{16}$. This algebra satisfies the conditions
of Proposition \ref{invform} and Theorem \ref{pareucl}. However, the
algebra and any non trivial subalgebra of it cannot be the algebra 
$\g_{\T}^*$ for a $3$-form $\T$ in sixteen variables. It turns out that
\bdm
\T(\spin(9) , \R^{16}) \, = \, 0 \, .
\edm
The proof is a longer algebraic computation and will be postponed  to the
appendix.
\end{exa}
\noindent
We remark that the results of this section cannot be generalized
directly to the case of $k$-forms. Spinorial connections related
with forms of higher degree behave differently.
Theorem \ref{pareucl} and Corollary \ref{parallelflach} are not true
for $4$-forms. Especially interesting is  dimension eight. 
A $4$-form $\T$ on $\R^8$ depends on $70$ parameters. On the other hand,
the $8$-dimensional spin representation is real and splits
$\Delta_8 = \Delta_8^+ \oplus \Delta_8^-$ into two $8$-dimensional
representations. Consider a spinor $\psi \in \Delta_8^+$ in one of these
components. The Clifford product $(X \haken \T) \cdot
\psi$ is a spinor in $\Delta_8^-$ and the condition
$\big(X \haken \T \big) \cdot \psi = 0$  for any  
$X \in \R^8$ 
is a system of $8 \cdot 8 = 64$ linear equations for the coefficients of
the $4$-form $\T$. Consequently, any spinor $\psi \in \Delta_8^{\pm}$
admits a family of $4$-forms $\T$ depending at least on $6$ parameters
such that $\mathfrak{g}_{\T}^* \cdot \psi = 0$. In fact, the number
of parameters is seven. Indeed, for any spinor
$\psi$, we consider the subspace
$\big\{ \T \in \Lambda^k(\R^n) \, : \ \mathfrak{g}_{\T}^* 
\cdot \psi = 0 \big\}$. 
It is invariant under the isotropy group of the spinor. In dimension eight, 
the isotropy group $\Spin(7)$ splits 
$\Lambda^4(\R^8)$ 
into four $\Spin(7)$-irreducible components of dimensions $1, \, 7, \, 27, 
\, 35$ (see \cite{Fern}).
In any case, there exist non trivial $4$-forms on $\R^8$ with non trivial 
parallel spinors. Since the
space $\Delta_8$ of all spinors in dimension eight coincides with the space
$\Delta_9$ of all spinors in dimension nine, we obtain $4$-forms in $\R^9$
with parallel spinors, too.
%
\section{$\nabla$-parallel $2$-forms on manifolds}
\noindent
Any metric connection on a Riemannian manifold defines several differential operators, like the Laplace
operator on forms or the Dirac operator on spinors. One can compare these
operators with the corresponding operator defined by the Levi-Civita 
connection. There is one particularly interesting formula
of that type, namely for the codifferential of an exterior form,
\bdm
\delta^{\nabla}\omega\ :=\ - \,\sum_{i=1}^n e_i \haken \nabla_{e_i}\omega\, .
\edm
We shall prove that the Riemannian
divergence of the torsion form coincides with its $\nabla$-divergence.
\begin{prop}\label{ext-diff}
Let $\nabla$ be a connection with skew-symmetric torsion. Then,
for any exterior form $\omega$, the following formula holds:
\bdm
\delta^{\nabla} \omega \ = \ \delta^g \omega \, - \, \frac{1}{2} \cdot
\sum_{i,j=1}^n (e_i \haken e_j \haken \T) \wedge (e_i \haken e_j \haken \omega)
\ .
\edm
In particular, for the torsion form itself, we obtain $\delta^{\nabla} \T = 
\delta^g \T$.
\end{prop}
\begin{proof} 
For simplicity, we prove the formula for $3$-forms. Then we get
\begin{eqnarray*}
\delta^{\nabla} \omega (X,Y) & = & - \, \sum_{i=1}^n \nabla_{e_i} 
\omega(e_i, X,Y) \ = \ - \, \sum_{i=1}^n e_i \big( \omega(e_i,X,Y)\big) \\
& & + \,  \sum_{i=1}^n \big[ \omega(\nabla_{e_i}e_i,X,Y) \, + \, \omega(e_i,
\nabla_{e_i}X,Y) \, + \, \omega(e_i,X,\nabla_{e_i}Y) \big]. \\
\end{eqnarray*}
Since the two connections are related by
$2 \cdot \nabla_XY - 2 \cdot \nabla^g_XY = \sum_{j=1}^n \T(X,Y,e_j)\cdot e_j$,
this can be rewritten in the form
\begin{eqnarray*}
\delta^{\nabla} \omega (X,Y)
& = & \delta^g \omega(X,Y) +  \frac{1}{2} \sum_{i,j=1}^n 
\big[\T(e_i,X,e_j)\, \omega(e_i,e_j,Y) +  \T(e_i,Y,e_j) \,
\omega(e_i,X,e_j) \big] \\
& = & \delta^g \omega(X,Y)  -  \frac{1}{2} 
\sum_{i,j=1}^n (e_i \haken e_j \haken \T) \wedge 
(e_i \haken e_j \haken \omega)(X,Y) \, . \quad\quad\quad\quad\quad\quad
 \qedhere
\end{eqnarray*}
\end{proof}
\begin{cor} 
If the torsion form $\T$ is $\nabla$-parallel, then its divergence vanishes, 
\bdm
\delta^g\T \ = \ \delta^{\nabla} \T \ = \ 0 \, .
\edm
\end{cor} 
\noindent
Let us discuss $\nabla$-parallel $2$-forms.
The differential equation reads as
\bdm 
\nabla^g_{\beta}\Omega_{\alpha\gamma} \ = \ \frac{1}{2}  
\sum^n_{j=1} \big\{ \Omega_{\gamma j} \cdot \T_{\beta j \alpha} \, - \, 
\Omega_{\alpha j} \cdot \T_{\beta j \gamma} \big\} \ .
\edm
Using the well known formulas for the exterior differential, the 
codifferential as well as for the Bochner-Laplace operator 
$\nabla^* \nabla$ we obtain
\begin{prop}
Let $\nabla$ be a metric connection $\nabla$
and skew-symmetric torsion. If $\Omega$ is a $\nabla$-parallel $2$-form, then
\begin{eqnarray*}
\delta^g \Omega & = & \frac{1}{2} \cdot \big(\Omega \haken \T \big) \ = 
\ \frac{1}{4}  \sum^n_{j,\beta, \gamma = 1} 
\Omega_{\beta j} \cdot \T_{\beta j \gamma} \cdot e_{\gamma} \, , \\
d \, \Omega & = & \sum^n_{j = 1} \big(e_j \haken \Omega\big) \wedge \big( e_j 
\haken \T) \\ 
& = & \frac{1}{6}  \sum^n_{j,\alpha, 
\beta, \gamma = 1} \big\{ \Omega_{\alpha j} \cdot \T_{\beta j \gamma} \, - \, 
\Omega_{\beta j} \cdot \T_{\alpha j \gamma} \, + \, \Omega_{\gamma j} 
\cdot \T_{\alpha j \beta} \big\} \cdot e_{\alpha} \wedge e_{\beta} \wedge
e_{\gamma} \, , \\
g\big( \Omega \, , \, \nabla^* \nabla^g \Omega \big) & = & 
\frac{1}{2}  \sum^n_{j,k,\alpha,\beta,\gamma = 1} \Omega_{\alpha \gamma}
\cdot \Omega_{\alpha k} \cdot \T_{\beta j k} \cdot \T_{\beta j \gamma} \ ,
\end{eqnarray*}
where $\nabla^* \nabla^g$ denotes the Riemannian Bochner-Laplace operator
acting on $2$-forms.
\end{prop}
\noindent 
In an adapted basis,  
$\Omega =  A_1 \cdot e_{1} \wedge e_2 + \cdots +  A_k \cdot
e_{2k-1} \wedge e_{2k}$,  
the third formula simplifies,
\bdm
g\big( \Omega \, , \, \nabla^* \nabla^g \Omega \big) \ = \ 
\frac{1}{2}  \sum_{\kappa = 1}^k \sum_{i, j = 1}^n 
\big( \T^2_{i j 2\kappa -1} \, + \, \T^2_{i j 2\kappa} \big) \cdot
A^2_{\kappa} \, .
\edm
It explains once again, from a geometric point of view, 
the proof of Proposition~\ref{invform}. We remark that there exist
indeed metric connections with skew-symmetric torsion and parallel
$2$-forms. Indeed, consider an almost hermitian manifold with totally 
skew-symmetric Nijenhuis tensor. Then there is a unique connection $\nabla$ 
preserving the hermitian structure with skew-symmetric torsion 
(see \cite{Friedrich&I1}). The fundamental
form of the hermitian structure is $\nabla$-parallel. A second example are
Sasakian manifolds. For these, the differential of the contact form is 
parallel with respect to the unique connection preserving the Sasakian 
structure.
%
\section{Schr\"odinger-Lichnerowicz type formulas for Dirac operators}
\label{fam-conn}\noindent
%
Consider a Riemannian spin manifold $(M^n,g, \mathrm{T})$ with $3$-form 
$\mathrm{T}$ as well as the one-parameter family of linear metric 
connections with torsion,
\bdm
\nabla^s_X Y\ :=\ \nabla^{g}_X Y + 2  s \cdot \mathrm{T}(X,Y,-)\,.
\edm
In particular, the superscript $s=0$ corresponds to the Levi-Civita 
connection, $\nabla^{g}\equiv\nabla^0$.
These connections can all be lifted to connections on the spinor
bundle $S$ of $M$, where they take the expression
\bdm
\nabla^s_X \psi\ :=\ \nabla^{g}_X \psi + s ( X\haken \mathrm{T}) \cdot 
\psi\,.
\edm
There is a formula for the square of the Dirac operator $D^s$ associated 
with the connection $\nabla^s$. In order to state it, let us introduce the
first order differential operator
\bdm
\D^s\psi \ :=\ \sum_{k=1}^n (e_k\haken\mathrm{T})\cdot \nabla^s_{e_k}\psi \ =\
\D^0\psi + s  \sum_k (e_k\haken \T)\cdot (e_k\haken \T)
\cdot\psi\,,
\edm
where $e_1,\ldots,e_n$ denotes an orthonormal basis. In fact, it will be 
convenient to use a separate notation for the algebraic $4$-form
derived from $\T$  appearing in the difference $\D^s-\D^0$:
\bdm
\sigma_{\T} \ :=\ \frac{1}{2} \sum_k (e_k\haken\T)\hut
(e_k\haken\T)  .
\edm
\begin{thm}[{\cite[Thm 3.1, 3.3]{Friedrich&I1}}]\label{FI-SL-AC}
Let $(M^n,g,\nabla^s)$ be an $n$-dimensional Riemannian
manifold with a metric connection $\nabla^s$ of skew-symmetric
torsion $4 \cdot s \cdot \T$. Then, the square of the Dirac operator
$D^s$ associated with  $\nabla^s$ acts on an arbitrary spinor field 
$\psi$ as
\be\label{FI-SL}
(D^s)^2\psi \ =\ \Delta^s(\psi) + 3  s \, d \T \cdot\psi
- 8  s^2\, \sigma_{\T}\cdot\psi+
2 s \, \delta \T \cdot\psi - 4  s \, \D^s\psi + \frac{1}{4}\,\Scal^s 
\cdot \psi,
\ee
where $\Delta^s$ is the spinor Laplacian of $\nabla^s$,
\bdm
\Delta^s(\psi)\ =\ (\nabla^s)^*\nabla^s\psi \ =\ -\sum_{k=1}^n \nabla^s_{e_k}
\nabla^s_{e_k} \psi +\nabla^s_{\nabla^{g}_{e_i} e_i}\psi\,.
\edm
Furthermore, the anticommutator of $D^s$ and $\omega$ is
\be\label{D-omega-anticomm}
D^s\cdot\T + \T\cdot D^s\ =\ d\T+\delta\T- 8 s \cdot \sigma_{\T}
-2 \, \D^s.
\ee
\end{thm}\noindent
$\Scal^s$ denotes the scalar curvature of the connection $\nabla^s$. Remark
that $\Scal^0 = \Scal^g$ is the usual scalar curvature of the underlying
Riemannian manifold $(M^n,g)$.

\smallskip\noindent
This formula for $(D^s)^2$ has the disadvantage of still
containing a first order differential operator as well as several
$4$-forms, which are difficult to treat algebraically. Inspired by
the homogeneous case,
we were looking for an alternative comparison of $(D^s)^2$ with the
Laplace operator of some \emph{other} connection $\nabla^{s'}$ from the same
family. For the computations, we need the square of $\T$ inside
the Clifford algebra. The proof of the following proposition is
completely similar to that of Proposition~3.1 in \cite{Agri} and will hence
be omitted.
\begin{prop}\label{3-form-square}
Let $\T$ be a $3$-form, and denote by the same symbol its
associated $(2,1)$-tensor. Then its square inside
the Clifford algebra has no contribution of degree $6$, and its scalar
and fourth degree part are given by
\bdm
\T^2_0\ =\ \frac{1}{6}\,\sum_{i,j=1}^n ||\T(e_i,e_j)||^2,\quad
\T^2_4\ =\ - \, 2 \cdot \sigma_{\T}. \qedhere
\edm
\end{prop}
\noindent
With these preparations in hand, we can state a more useful
Schr\"odinger-Lichnerowicz type formula for $(D^s)^2$. It links the
Dirac operator for the parameter $s/3$ with the Laplacian for the
parameter $s$. The remainder is a zero order
operator. Similar formulas can be found in \cite{Bis} and, for  homogeneous
spaces, in \cite{Agri}.
\begin{thm}\label{new-weitzenboeck}
The spinor Laplacian $\Delta^s$ and the square of the Dirac operator
$D^{s/3}$ are related by
\bdm
(D^{s/3})^2\ =\ \Delta^s + s \cdot d\T 
+\frac{1}{4}\,\Scal^g - 2s^2 \cdot \T^2_0.
\edm
\end{thm}
\begin{proof}
By the formula from Theorem~\ref{FI-SL-AC},
\bdm
(D^s)^2 + 4s \, \D^s \ =\ \Delta^s + 3s \, d\T 
- 8s^2 \, \sigma_{\T}   +
 2s \, \delta \T   + \frac{1}{4}\,\Scal^s.
\edm
But since $D^s=D^0 +3s \cdot \T$, the left hand side can equally be
rewritten
\bdm
(D^s)^2 +  4s \, \D^s \ =\ (D^0)^2+ 3s  (\T D^0 + D^0\T) +
9s^2 \, \T^2   + 4s \, \D^s .
\edm
We use equation~(\ref{D-omega-anticomm}) to express $\D^s$ by
the anticommutator $\T D^s + D^s\T$:
\bdm
2 \, \D^s \ = \ d\T + \delta\T - 8s \cdot \sigma_{\T} - (D^0\T + \T D^0) 
- 6s \cdot \T^2 \ .
\edm 
Now we obtain
\begin{eqnarray*}
(D^s)^2  + 4s \, \D^s & =& (D^0)^2 + s (\T D^0+D^0\T)-
 3s^2 \, \T^2 -16s^2\cdot\sigma_{\T}+
2s \cdot d\T \\
& = & (D^{s/3})^2 - 4s^2 \cdot \T^2 - 16s^2\cdot\sigma_{\T}+
2s \cdot d\T \ .
\end{eqnarray*}
We observe that $D^{s/3}$ hence appears by quadratic completion.
Now it suffices to insert this result in the formula of Theorem 
\ref{FI-SL-AC} and to use Proposition~\ref{3-form-square} as well
as well as the easy relation between scalar curvatures, 
$\Scal^s = \Scal^g - 24 s^2 \, \T^2_0$.
\end{proof}
\noindent
Integrating the latter formula on a compact manifold $M^n$, we obtain
\bdm
\int_{M^n}  ||D^{s/3} \psi||^2 \ = \
\int_{M^n} \Big[ ||\nabla^{s} \psi||^2 + s \langle d\T \cdot \psi, 
\psi \rangle + \frac{1}{4} \Scal^g \cdot ||\psi||^2 - 2s^2 \, \T^2_0 \cdot ||\psi||^2 \Big] \, . 
\edm\

\noindent
A first consequence is a non linear version of Corollary \ref{parallelflach}.
\begin{thm}\label{parallelflach2}  
Let $(M^n, g, \T)$ be a compact, Riemannian spin
manifold of non positive scalar curvature, $\Scal^g \leq 0$, and suppose that 
the $4$-form $d\T$ acts on spinors as a non positive endomorphism. If there
exists a solution $\psi \neq 0$  of the equation
\bdm
\nabla^{\T}_X\psi \ = \ \nabla^g_X \psi \, + \, (X \haken \T) \cdot \psi \ = \ 0 \ ,
\edm
then the $3$-form and the scalar curvature vanish, $\T = 0 = \Scal^g$, 
and $\psi$ is parallel with respect to the Levi-Civita connection.
\end{thm}
\begin{NB} 
Let us compare  Theorem \ref{parallelflach2} with 
the integral formula in \cite{Friedrich&I1}. There, we need the 
condition that $d\T + 8 \cdot \sigma_{\T}$ is a non positive endomorphism
in order to prove the same result. Since $\sigma_{\T}$ is neither positive 
nor negative, the two conditions
are independent. The advantage of Theorem \ref{parallelflach2} is that
only the algebraic type of the exterior differential $d \T$ is involved,
but not the algebraic type of the torsion form $\T$ itself (see the proof
of Theorem \ref{pareucl}).
\end{NB}
\noindent
Theorem \ref{parallelflach2} applies, in particular, to Calabi-Yau or 
Joyce manifolds. These
are compact, Ricci-flat Riemannian manifolds in dimensions $n=6,7$ with one 
parallel spinor field.
Let us perturb the connection $\nabla^g$ by a $3$-form such that $d\T$ is non 
positive on spinors. Then the
new connection $\nabla^{\T}$ does not admit $\nabla^{\T}$-parallel spinor 
fields. Nilmanifolds and their compact quotients 
$M^n = \mathrm{G}/\Gamma$ are a
second family of examples where the theorem applies. A further family of
examples arises from certain naturally reductive spaces and a torsion form
$\T$ being proportional to the torsion form of the canonical connection, see
\cite{Agri}.
%
\section{$1$-parameter families of connections with parallel spinors}
%
\noindent
Consider a triple $(M^n,g,\T)$ consisting of a Riemannian manifold together
with a fixed $3$-form $\T \neq 0$. Let us ask for parameters $s_0$ such that 
the connection $\nabla^{s_0}$ admits a parallel spinor. The first example
describes a case with parallel spinors for more then only one parameter
in the family. 
\begin{exa}\label{Lie-grp} 
Let $G$ be a simply connected Lie group, $g$ a biinvariant metric and consider
the torsion form $\T(X,Y,Z) := g([X,Y],Z)$. The connections $\nabla^{\pm 1/4}$
are flat (see \cite{KN2}). In particular, there are $\nabla^{\pm 1/4}$-parallel
spinor fields.
\end{exa}
\noindent
The integrability condition of Theorem
\ref{FI-SL-AC} implies that the function 
\bdm
\mathrm{G}(m, s) \ := \ \mathrm{Det} \big\{3  s \, d \T 
- 8  s^2\, \sigma_{\T} + 2 s \, \delta \T  + \frac{1}{4}\,\Scal^s \big\}(m)
\edm
vanishes at $s_0$ and all $m \in M^n$. Here we treat  forms 
as endomorphisms acting on spinors. The function $\mathrm{G}(s)$ is a 
polynomial. If the Riemannian scalar curvature is not identically zero,  
there exists only a finite number of parameters with $\nabla^s$-parallel 
spinors. Let us discuss low dimensions.
\begin{exa}[The $3$-dimensional case]
Consider the $3$-dimensional sphere $(S^3, g, dS^3)$ equipped
with its standard metric and the volume form $\T = dS^3$. The equation
\bdm
\nabla_X^s \psi \ = \ \nabla_X^g \psi \, +  \, s \cdot (X \haken \T) \cdot 
\psi \ = \ \nabla_X^g \psi \, + \, s \cdot X \cdot \psi \ = \ 0
\edm
is the usual Killing spinor equation. There are solutions on the 
$3$-dimensional sphere for both parameters $ s = \pm 1/2$.
In dimension $n=3$, this is the only manifold admitting parallel spinors 
with respect to a  non trivial
$3$-form. Indeed, any $\T$ is proportional to the volume form, 
$\T = f \cdot dM^3$, where $f$ is a real-valued smooth function on $M^3$.
If the equation
\bdm
\nabla_X^g \psi \, + \, (X \haken \T)\cdot \psi \ = \ 
\nabla_X^g \psi \, + \, f \cdot X \cdot \psi \ = 0
\edm
admits a non trivial solution $\psi$, then by a Theorem of A.~Lichnerowicz
(see \cite{Li}) $f$ is constant and $(M^3,g)$ is a space form. 
\end{exa}
\noindent
In dimension four, we split any $2$-form $\omega \in
\Lambda^2(\R^4)$ into its self-dual and 
anti-self-dual part, $\omega = \omega_+ + \omega_-$.
\begin{lem} 
An element $a + \omega + f \cdot e_{1234} \in \mathrm{Cl}(\R^4)$
acts on the space $\Delta_4$ of spinors and its determinant is given by the
formula
\bdm
\mathrm{Det}\big[a + \omega + f \cdot e_{1234}\big] \ = \
\big[(a + f)^2 \, + \, 4 \cdot ||\omega_+||^2 \big]\cdot
\big[(a - f)^2 \, + \, 4 \cdot ||\omega_-||^2 \big] \ .
\edm
For any $3$-form $\T \in \Lambda^3(\R^4)$ the corresponding $4$-form
$\sigma_{\T}$ vanishes, $\sigma_{\T} = 0$.
\end{lem}  
\begin{proof} 
Any $3$-form in $\R^4$ is $\SO(4)$-equivalent to the form
$a \cdot e_1 \wedge e_2 \wedge e_3$ containing only one summand. This
normal form implies $\sigma_{\T} = 0$ immediately. The formula for the
determinant follows from a matrix representation of the Clifford algebra. 
\end{proof}
\noindent
The equation $\mathrm{G}(m,s) = 0$ yields the following
condition {\it not} expressing the
full integrability conditions for the existence of a parallel spinor. 
\begin{prop} Let $(M^4,g,\T)$ be a Riemannian spin manifold equipped with
a $3$-form $\T$. If the connection $\nabla^s$ admits a non trivial parallel 
spinor,  the following equations hold at any point:
\begin{enumerate}
\item $ 12 \cdot s \cdot d \T \ = \ \pm \, ( \Scal^g \, - \, 24 \cdot s^2 \cdot \T_0^2) \cdot 
dM^4$.
\item $\delta(\T)$ is a (anti)-self-dual $2$-form.
\end{enumerate}
\end{prop}
\begin{exa} 
Using the unique $3$-dimensional example $S^3$ and its Killing
spinors, we obtain by $M^4 := S^3 \times \R^1$ and $\T := dS^3$ an example in dimension four. Indeed, the $3$-dimensional Killing
spinors are $\nabla^{\pm 1/2}$-parallel on $M^4$. They do not depend on the 
$\R^1$-coordinate. 
\end{exa}
\noindent
The integrability condition restricts the admissible
parameters via a polynomial equation involving the scalar curvature
and the torsion form of the triple $(M^n,g,\T)$. Globally, not all of these
values are possible.
\begin{thm}\label{parameter}  
Let $(M^n,g,\T)$ be a compact triple. For any $\nabla^s$-parallel spinor 
$\psi$, the following formula holds: 
\bdm
64 \cdot s^2 \int_{M^n} \langle \sigma_{\T} \cdot \psi \, , \, \psi 
\rangle \, + \, \int_{M^n} \mathrm{Scal}^s \ = \ 0 \ .
\edm
If the mean value of $\langle \sigma_{\T}  \cdot \psi \, , \, \psi \rangle$ 
does not vanish,  the parameter $s$ is given by 
\bdm
s \ = \ \frac{1}{8} \int_{M^n} \langle d\T \cdot \psi \, , \, \psi 
\rangle \Big/ 
\int_{M^n} \langle \sigma_{\T} \cdot \psi \, , \, \psi \rangle  
\edm
If the mean value of $\langle \sigma_{\T}  \cdot \psi \, , \, \psi \rangle$ 
vanishes,  the parameter $s$ depends only on the Riemannian scalar
curvature and on the length of the torsion form, 
\bdm
0 \ = \ \int_{M^n} \mathrm{Scal}^s \ = \ \int_{M^n} \Scal^g \, - \,
24s^2 \int_{M^n} \T^2_0 \ .
\edm
Finally, if the $4$-forms $d\T$ and $\sigma_{\T}$ are proportional,
there are at most three parameters with $\nabla^s$-parallel spinors.   
\end{thm}
\begin{proof} 
We use the integrability conditions for parallel spinors
from Theorem \ref{FI-SL-AC}. Let $\psi$ be a $\nabla^s$-parallel spinor
of length one. Then we obtain
\bdm
3s \int_{M^n} \langle d\T \cdot \psi\, , \, \psi \rangle \, - \, 
8s^2 \int_{M^n} \langle \sigma_{\T} \cdot \psi\, ,\, \psi \rangle \, + \, 
\frac{1}{4} \int_{M^n} \Scal^g \, - \, 6s^2 \int_{M^n} \T^2_0 \ = \ 0 \ . 
\edm
On the other side, the anticommutator relation between $D^s$ and $\T$
as well as the symmetry property of the Dirac operator in $\mathrm{L}^2$
yields
\bdm
0 \ = \ \int_{M^n} \langle \T \cdot \psi \, , \, D^s \psi \rangle \ = 
\ \int_{M^n} \langle D^s \T \cdot \psi \, , \, \psi \rangle \ = \ \int_{M^n} 
\langle d\T \cdot \psi \, , \, \psi \rangle \, - \,  8s \int_{M^n} 
\langle \sigma_{\T} \cdot \psi \, , \, \psi \rangle \ .
\edm
If the mean values of $\langle \sigma_{\T} \cdot \psi \, , \, \psi \rangle$ 
does not
vanish, then the second equation determines the parameter $s$,
\bdm
s \ = \ \frac{1}{8} \int_{M^n} \langle d\T \cdot \psi \, , \, \psi \rangle 
\Big/ \int_{M^n} \langle \sigma_{\T} \cdot \psi \, , \, \psi \rangle \ .
\edm
If the mean values of $\langle \sigma_{\T} \cdot \psi \, , \, \psi \rangle$ 
vanishes,
then the mean value of $\langle d \T \cdot \psi \, , \, \psi \rangle$ 
vanishes, too.
The first formula yields the result.
\end{proof}
\begin{NB} 
In Proposition~\ref{TandminusT}, we discuss an example of a non-flat
connection on the compact, $7$-dimensional Aloff-Wallach space $N(1,1)$ 
such that $\nabla^{s_0}$ \emph{and} $\nabla^{-s_0}$ admit parallel spinors
for suitable $s_0$, 
hence showing that both cases from Theorem~\ref{parameter} can actually
occur in non-trivial situations. The "trivial" cases we knew about
before are, of course, Lie groups (Example~\ref{Lie-grp}). 
Example \ref{scal-null} illustrates
how a parallel spinor can occur for zero scalar curvature and $d \T$ 
proportional to $\sigma_{\T}$. In the same vein, we construct on 
$N(1,1)$ a spinorial connection defined by a $4$-form $R$ such that
$\nabla^{R}$  and $\nabla^{-R}$ admit parallel spinors 
(Proposition~\ref{RandminusR}). 
\end{NB}
\noindent
If the torsion form $\T$ of the linear connection
is $\nabla$-parallel, we have 
$d \T = 2 \cdot \sigma_{\T}$ and $\delta(\T) = 0$.
This situation occurs if $M^n = G/H$ 
is a reductive space and $\T$ is the torsion form of its natural
connection (see \cite{KN}) or for Sasakian manifolds, nearly K\"ahler
manifolds, nearly parallel $\mathrm{G}_2$-manifolds equipped with their
unique connection preserving the corresponding geometric structure
(see \cite{Friedrich&I1}). A direct consequence of Theorem 
\ref{new-weitzenboeck} and Theorem 
\ref{parameter} is the following
\begin{cor} 
Let $(M^n,g,\T)$ be a compact Riemannian manifold and suppose that
the exterior differential of the $3$-form $\T$ is proportional to the
$4$-form $\sigma_{\T}$, $d\T = 2 \cdot \sigma_{\T}$. If a connection
$\nabla^s$ with $s \neq 1/4$ admits a parallel spinor field,  the 
first eigenvalue of the Dirac operator $D^{s/3}$ is bounded by
\bdm
6 \cdot \mathrm{vol}(M^n,g) \cdot \lambda^2_1(D^{s/3}) \ \leq \ \int_{M^n} 
\Scal^g \ .
\edm
If $\nabla^{1/4} = \nabla$ admits a parallel spinor field, 
\bdm
\mathrm{vol}(M^n,g) \cdot \lambda^2_1(D^{1/12}) \ \leq \ \frac{1}{8}
\int_{M^n} \Scal^g \, + \, \frac{1}{16} \int_{M^n} \T_0^2 \ .
\edm
\end{cor}
\begin{NB}
On a naturally reductive space $M=G/K$, $\nabla^{1/4}$ is the canonical
connection and \cite[Corollary 3.1.]{Agri} shows that a  
$\nabla^{1/4}$-parallel spinor realizes indeed this lower bound for 
$\lambda^2_1(D^{1/12})$ provided the Casimir operator $\Omega_{\g}$
is non-negative. 
\end{NB}
\noindent
In case that the torsion form of the triple $(M^n,g,\T)$ arises 
from some special 
non-integrable geometric structure (see \cite{Fri2}), then usually 
only one connection in the family admits $\nabla^s$-parallel spinors.
A uniqueness of that type requires additional arguments involving
the special geometric structure. For example, consider a 
$5$-dimensional Sasakian manifold
$(M^5, g, \xi, \eta , \varphi)$. Denote by $\nabla$ its unique connection
with skew-symmetric torsion and preserving the contact structure.
Its torsion is given by the formula $\T = \eta \wedge d \eta$ (see \cite{Friedrich&I1}). In an adapted local frame, we have the formulas
\bdm
\T \ = \ \eta \wedge d \eta \ = \ 2 \cdot (e_{12} + e_{34}) \wedge e_5 , 
\quad d \eta \ = \ 2 \cdot (e_{12} + e_{34}) .
\edm 
We consider the family $\nabla^s$
of connections. The first admissible case $s = 1/4$ is the connection
$\nabla$ preserving the contact structure we started with. In the papers 
\cite{Friedrich&I1}, \cite{Friedrich&I2} the integrability 
conditions for $\nabla^{1/4}$-parallel spinors have been discussed completely.
In particular, there are compact examples. 
For any Sasakian structure, we have
\bdm
2 \cdot \nabla^g_X \xi \ = \  X \haken d\eta \, . 
\edm
Suppose that there exists a $\nabla^s$-parallel spinor $\psi^*$ for some 
parameter $s \neq 1/4$.
We introduce the vector field $\xi^*$ via the algebraic equation 
$\xi^* \cdot \psi^* = i \cdot \psi^*$. Then $\xi^*$ is $\nabla^s$-parallel,
\bdm
\nabla^g_X \xi^* \ = \ - \, 2s \cdot \T(X,\xi^*, -) \, .
\edm 
Let us consider the inner product $f := g (\xi \, , \,  \xi^* )$ of 
the two vector fields. Its differential is given by the formula
\bdm
2 \cdot df \ = \ (4s \, - \, 1) \, \xi^* \haken d \eta \, .
\edm
In particular, $\xi(f) = \xi^*(f) = 0$. Next, we compute the commutator
of the vector fields
\bdm
\big[\xi \, , \, \xi^* \big] \ = \ \nabla^g_{\xi}\xi^* \, - \, 
\nabla^g_{\xi^*}\xi \ = \ - \, (1/2 \, + \, 2s) \, \xi^* \haken d \eta \  
= \ - \, \frac{4s + 1}{4s - 1} \cdot df \, .
\edm
First we discuss the case that $s \neq \pm 1/4$. Since $\big[\xi, \xi^*\big](f) = 0$, we conclude that $||\mathrm{grad}(f)||^2
= 0$ holds and then we obtain $\xi^* \haken d \eta = 0$.
Consequently, $\xi^*$ is proportional to the vector field $\xi$. In particular,
 $\xi$ is $\nabla^s$-parallel,
\bdm
\nabla^g_X \xi \ = \ -\, 2s \cdot \T(X,\xi, -) \ = \ 2s \cdot X \haken 
d \eta \, .
\edm
If $s \neq 1/4$, the latter equation contradicts the differential equation
for the Killing vector field $\xi$ of a Sasakian structure. Finally, we
study the remaining case  $s = - \, 1/4$. Then we have
\bdm
3s \cdot d \T \, - \, 8s^2 \cdot \sigma_{\T} \, + \, \frac{1}{4} 
\cdot \Scal^s \ 
= \ - \, 8 \cdot e_{1234} \, + \ \frac{1}{4} \cdot \Scal^{- 1/4} .
\edm 
The endomorphism $e_{1234}$ acts on spinors with constant eigenvalues $\pm 1$.
Therefore, if $\psi^*$ is a $\nabla^{- 1/4}$-parallel spinor,  the
scalar curvature $\Scal^{- 1/4}$ is constant and $\psi^*$ is an eigenspinor
of this endomorphism, $e_{1234} \cdot \psi^* = \epsilon \cdot \psi^*$. 
Since the connection $\nabla^{1/4}$ preserves the contact structure,
the covariant derivative $\nabla_X^{1/4} \psi^*$ satisfies the same algebraic
equation. With respect to
\bdm
0 \ = \ \nabla_X^{- 1/4} \psi^* \ = \ \nabla_X^{1/4} \psi^* \, - \, 
\frac{1}{2} \cdot (X \haken \T) \cdot \psi^*
\edm
we conclude that for any vector $X$ the spinor $\psi^*$ satisfies the
equation
\bdm
e_{1234} \cdot (X \haken \T) \cdot \psi^* \ = \ \epsilon \cdot (X \haken \T) 
\cdot \psi^* \, .
\edm
Inserting $X = e_1$ we obtain
$e_{1234} \cdot e_{25} \cdot \psi^* =  \epsilon \cdot e_{25} \cdot \psi^*$  
and $e_{1234} \cdot \psi^* = \epsilon \cdot \psi^*$.
The relations in the Clifford algebra yield immediately that $\psi^* = 0$.
All together, we proved:
\begin{prop} Let $(M^5, g, \xi, \eta , \varphi)$ be a $5$-dimensional
Sasakian manifold and denote by $\nabla$ its unique connection
with skew-symmetric torsion $\T$ and preserving the contact structure. 
If a connection $\nabla^s$ in the family through $\nabla$ 
admits a parallel spinor field, then $s = 1/4$ and the connection
is $\nabla$. 
\end{prop}
%
\section{Torsion forms with parallel spinors on Aloff-Wallach spaces}
\label{AW-exa}\noindent
%
The goal of this section is to construct on the Aloff-Wallach manifold 
$N(1,1) =\SU(3)/S^1$  a two-parameter family of metrics $g=g_{s,y}$ 
that admits, for every
$g_{s,y}$, two inequivalent cocalibrated $\mathrm{G}_2$-structures. Moreover, we investigate the torsion forms of their unique connections (see \cite{Friedrich&I1}) 
as well as other geometric data of these connections.
We use the computations available in \cite[p.109 ff]{BFGK}, 
which we hence shall not reproduce here. Consider the embedding
$S^1\ra \SU(3)$ given by $e^{i\theta}\mapsto \diag(e^{i\theta}, e^{i\theta},
e^{-2i\theta})$. The Lie algebra $\su(3)$ splits into $\su(3)=\m+\R$,
where $\R$ denotes the Lie algebra of $S^1$ deduced from the given embedding.
The space $\m$ has a preferred direction, namely the subspace
$\m_0$ generated by the matric $L:=\diag(3i, -3i,0)$. Let $E_{ij}\, (i<j)$ 
be the matrix with $1$ at the place $(i,j)$ and zero 
elsewhere, and define $A_{ij}=E_{ij}-E_{ji},\,\tilde{A}_{ij}=i(E_{ij}+E_{ji})$.
We set $\m_1:=\mathrm{Lin}\{A_{12}, \tilde{A}_{12}\}$,  
$\m_2:=\mathrm{Lin}\{A_{13}, \tilde{A}_{13}\}$ and 
$\m_3:=\mathrm{Lin}\{A_{23}, \tilde{A}_{23}\}$. The sum
$\m_1\op\m_2\op\m_3$ is an algebraic
complement of $\m_0$ inside $\m$, and in fact all spaces $\m_i$ are
pairwise perpendicular  with respect to the Killing form 
$B(X,Y):=-\mathrm{Re}(\tr XY)/2$. Hence, the following formula
\bdm
\g_{s,y} \ :=\ \frac{1}{s^2}\, B\big|_{\m_0}\, + \, B\big|_{\m_1}\,+\,
 \frac{1}{y}B\big|_{\m_2}\,+\,\frac{1}{y}B\big|_{\m_3}
\edm
defines a two-parameter family of metrics on $N(1,1):=\SU(3)/S^1$. It 
is a subfamily
of the family considered in \cite[p.109 ff]{BFGK}; in particular, 
($s=1,y=2$) corresponds to the $3$-Sasakian metric that has three
Killing spinors with Killing number $1/2$, and ($s=1,y=2/5$) is the 
Einstein metric with one Killing spinor with Killing number $-3/10$
(see \cite[Thm 12, p.116]{BFGK}). An orthonormal basis of $\m$ is given by
\bdm
X_1\, =\,A_{12},\ X_2\, =\, \tilde{A}_{12},\
X_3\, =\,\sqrt{y} A_{13},\ X_4\, =\, \sqrt{y}\tilde{A}_{13},\
X_5\, =\,\sqrt{y}A_{23},\ X_6\, =\, \sqrt{y}\tilde{A}_{23},
\edm
and $X_7 = s\cdot L/3$.
The isotropy representation $\Ad(\theta)$ leaves the vectors
$X_1,X_2$ and $X_7$ invariant, and acts as a rotation by
$3\theta$ in the $(X_3,X_4)$-plane and in the $(X_5,X_6)$-plane. 
We use the standard realization of the
$8$-dimensional $\Spin(7)$-representation $\Delta_7$ as given in 
\cite[p.97]{BFGK} or \cite[p.13]{Dirac-Buch}, and denote by
$\psi_i,i=1,\ldots 8$ its basis ($u_i$ in the notation of \cite{BFGK}). 
One then checks that $\psi_3,\psi_4,\psi_5$ and $\psi_6$
are fixed under the lift $\tilde{\Ad}(\theta)$ of the isotropy representation 
to $\Spin(7)$. Thus, they define constant sections in the spinor bundle
$S=\SU(3)\x_{\tilde{\Ad}}\Delta_7$.  The Levi-Civita connection
of $N(1,1)$ is described by a map $\Lambda:\ \m\mapsto\so(7)$, whose lift
$\tilde{\Lambda}:\ \m\mapsto\spin(7)$  is given by\footnote{Notice the 
following typo in the reference: the right definition of $d$ on page 112
is $d:=\sqrt{xy/z}+\sqrt{yz/x}-\sqrt{xz/y}$. For our metrics, $x=1$ and 
$y=z$.} (\cite[p.112]{BFGK})
\begin{eqnarray*}
\tilde{\Lambda}(X_1) &= & +\,\frac{1}{2s}\,e_2\cdot e_7 - 
\left[\frac{1}{2}-\frac{y}{4}\right][e_3\cdot e_5+ e_4\cdot e_6]\\ 
\tilde{\Lambda}(X_2) &=& -\,\frac{1}{2s}\,e_1\cdot e_7 - 
\left[\frac{1}{2}-\frac{y}{4}\right][e_4\cdot e_5-e_3\cdot e_6]\\
\tilde{\Lambda}(X_3) &= &+ \,\frac{y}{4s}\,e_4\cdot e_7 - \frac{y}{4}\,
[e_2\cdot e_6- e_1\cdot e_5]\\ 
\tilde{\Lambda}(X_4) &= &-\,\frac{y}{4s}\,e_3\cdot e_7 + 
\frac{y}{4}[e_1\cdot e_6+e_2\cdot e_5] \\  
\tilde{\Lambda}(X_5) &= &-\,\frac{y}{4s}\,e_6\cdot e_7 - 
\frac{y}{4}[e_1\cdot e_3+ e_2\cdot e_4]\\ 
\tilde{\Lambda}(X_6) &=& +\,\frac{y}{4s}\,e_5\cdot e_7 - 
\frac{y}{4}[e_1\cdot e_4- e_2\cdot e_3] \\
\tilde{\Lambda}(X_7) &= & \,\frac{s}{2}[2 e_1\cdot e_2 + e_3\cdot e_4 - 
e_5\cdot e_6] -\frac{1}{2s}e_1\cdot e_2-\frac{y}{4s}e_3\cdot e_4 +
\frac{y}{4s}e_5\cdot e_6.
\end{eqnarray*}
We now make the following Ansatz for an algebraic $3$-form on $\m$,
\begin{eqnarray*}
\T & = & \ \alpha\, X_1\wedge X_3\wedge X_5\, +\,
\beta\, X_1\wedge X_4\wedge X_6\,+\,\gamma X_2\wedge X_4\wedge X_5\, +\,
\delta\,X_2\wedge X_3\wedge X_6\\
&+ &\mu\, X_1\wedge X_2\wedge X_7\,+\,
\nu\, X_3\wedge X_4\wedge X_7\,+\, \eta\, X_5\wedge X_6\wedge X_7. 
\end{eqnarray*}
For notational convenience, we shall
write $X_{ijk}$ for $X_i\wedge X_j\wedge X_k$, and similarly for forms of
any degree.
In order to define a global form on $N(1,1)$, an algebraic  form on $\m$ needs
to be invariant under the isotropy representation. This is true for 
$X_{127}$, $X_{347}$, and $X_{567}$, whereas for example $X_{135}$ does not
exist globally. However, one easily checks that  the two $2$-forms
$X_{35}+X_{46},\ X_{45}-X_{36}$
are isotropy invariant, and this will suffice to check that
all forms to follow are indeed well-defined on $N(1,1)$.
In any event, $X_1\haken\T$ acts on algebraic spinors by Clifford 
multiplication
with $\alpha\, e_3\cdot e_5+\beta\, e_4\cdot e_6+\mu\, e_2\cdot e_7$, and
similarly for $X_2,\ldots,X_7$. 
\begin{prop}
The spinor field $\psi_3$ satisfies the equation 
$\nabla^g_X \psi_3+(X\haken\T)\cdot \psi_3=0$ exactly for one $3$-form
$\T:=\T_3$,  
\begin{eqnarray*}
\T_3 & := & \left[\frac{1}{2}-\frac{y}{4}+\frac{1+y}{6s}-\frac{s}{3} \right][X_{135}\, +\, X_{146}] +
\left[\frac{1}{2}- \frac{y}{4}-\frac{1+y}{6s}+\frac{s}{3} \right][X_{245}\, -\,X_{236}] \\
&+& \left[\frac{2y-1}{6s}-\frac{2s}{3}\right] \,X_{127}
+ \left[\frac{4+y}{12s}-\frac{2s}{3}\right][X_{347}\,-\,X_{567}].
\end{eqnarray*}
\end{prop}
\begin{proof}
A computer computation yields that the overdetermined system of 
equations $\nabla^g_{X_i} \psi_3+(X_i\haken\T)\cdot \psi_3=0$ 
(for $i=1,\ldots,7$)  reduces to a linear system of seven equations in the
seven variables $\alpha,\ldots,\eta$ with two free parameters $s,y>0$:
\bdm
\frac{1}{2s}+1-\frac{y}{2}-\alpha-\beta+\mu=0,\ -\frac{1}{2s}+1-\frac{y}{2}-\gamma+\delta-\mu=0,\ 
\frac{y}{4s}-\alpha-\delta+\nu=0,\ \frac{y}{4s}-\beta+\gamma+\nu=0,\ 
\edm
\bdm
-\frac{y}{4s}+\alpha-\gamma+\eta=0,\ -\frac{y}{4s}+\beta+\delta+\eta=0,\
2s-\frac{1+y}{2s}+\mu+\nu-\eta=0. 
\edm
One then verifies that the coefficients given in the proposition
are its unique solution.
\end{proof}
\begin{prop}
The spinor field $\psi_4$ satisfies the equation 
$\nabla^g_X \psi_4+(X\haken\T)\cdot \psi_4=0$ exactly for one $3$-form 
$\T:=\T_4$, 
\begin{eqnarray*}
\T_4 & := & \left[\frac{1}{2}-\frac{y}{4}-\frac{1+y}{6s}+\frac{s}{3} \right][X_{135}\, +\, X_{146}] +
\left[\frac{1}{2}-\frac{y}{4}+ \frac{1+y}{6s}-\frac{s}{3} \right][X_{245}\, -\,X_{236}]\\
& + & \left[\frac{2y-1}{6s}-\frac{2s}{3}\right] \,X_{127}
+ \left[\frac{4+y}{12s}-\frac{2s}{3}\right][X_{347}\,-\,X_{567}].
\end{eqnarray*}
\end{prop}
\begin{proof}
The linear system determined by
$\nabla^g_{X_i} \psi_4+(X_i\haken\T)\cdot \psi_4=0$ reads as
\bdm
\frac{1}{2s}-1+\frac{y}{2}+\alpha+\beta+\mu=0,\ 
\frac{1}{2s}+1-\frac{y}{2}-\gamma+\delta+\mu=0,\
\frac{y}{4s}+\alpha+\delta+\nu=0,\ 
\frac{y}{4s}+\beta-\gamma+\nu=0,
\edm
\bdm
\frac{y}{4s}+\alpha-\gamma-\eta=0,\ 
\frac{y}{4s}+\beta+\delta-\eta=0,\ 
-2s+\frac{1+y}{2s}-\mu-\nu+\eta=0 .
\edm
Its unique solution leads to the formulas above.
\end{proof}
\begin{prop}
The spinor field $\psi_5$ satisfies the equation 
$\nabla^g_X \psi_5+(X\haken\T)\cdot \psi_5=0$ exactly for one $3$-form 
$\T:=\T_5$, 
\begin{eqnarray*}
\T_5 & := & \left[\frac{1}{6}+\frac{y}{12}+\frac{y-1}{6s}\right] 
\,[X_{135}\, +\, X_{146}\, +\,X_{245}\, -\,X_{236}]\, +\,
\left[\frac{2}{3}- \frac{2y}{3}-\frac{2y+1}{6s}\right]\,X_{127}\\
& + & \left[ \frac{1}{3}-\frac{y}{3}-\frac{4-y}{12s}\right]\, [ X_{347} \, 
- \, X_{567}].
\end{eqnarray*}
\end{prop}
\begin{proof}
The linear system $\nabla^g_{X_i} \psi_5+(X_i\haken\T)\cdot \psi_5=0$ is
of slightly different type,
\bdm
\frac{y}{4s}+\frac{y}{2}-\alpha+\delta+\nu=0,\
\frac{y}{4s}+\frac{y}{2}-\beta-\gamma+\nu=0,\ 
\frac{y}{4s}+\frac{y}{2}-\alpha-\gamma-\eta=0,\ 
\frac{y}{4s}+\frac{y}{2}-\beta+\delta-\eta=0,\ 
\edm
\bdm
\frac{1}{2s}-1+\frac{y}{2}+\gamma-\delta+\mu=0,\ 
\frac{1}{2s}-1+\frac{y}{2}+\alpha+\beta+\mu=0,
\frac{y-1}{2s}+ \mu-\nu+\eta=0.
\edm
The main reason for this is that $\psi_5$ and $\psi_6$ span the kernel of 
the  first summand of $\tilde{\Lambda}(X_7)$, hence the last equation
contains no term linear in $s$.
\end{proof}
\begin{prop}
The spinor field $\psi_6$ satisfies the equation 
$\nabla^g_X \psi_6+(X\haken\T)\cdot \psi_6=0$ exactly for one $3$-form 
$\T:=\T_6$,
\begin{eqnarray*}
\T_6 & := & \left[\frac{1}{6}+\frac{y}{12}-\frac{y-1}{6s}\right] 
\,[X_{135}\, +\, X_{146}\, +\,X_{245}\, -\,X_{236}]\, +\,
\left[-\frac{2}{3}+ \frac{2y}{3}-\frac{2y+1}{6s}\right]\,X_{127}\\
& + & \left[- \frac{1}{3}+\frac{y}{3}-\frac{4-y}{12s}\right]\, [ X_{347}
\, - \,  X_{567}].
\end{eqnarray*}
\end{prop}
\begin{proof}
The linear system $\nabla^g_{X_i} \psi_6+(X_i\haken\T)\cdot \psi_6=0$ is
\bdm
\frac{y}{4s}-\frac{y}{2}+\alpha-\delta+\nu=0,\
\frac{y}{4s}-\frac{y}{2}+\beta+\gamma+\nu=0,\ 
\frac{y}{4s}-\frac{y}{2}+\alpha+\gamma-\eta=0,\ 
\frac{y}{4s}-\frac{y}{2}+\beta-\delta-\eta=0,\ 
\edm
\bdm
-\frac{1}{2s}-1+\frac{y}{2}+\gamma-\delta-\mu=0,\ 
-\frac{1}{2s}-1+\frac{y}{2}+\alpha+\beta-\mu=0,
\frac{1-y}{2s}- \mu+\nu-\eta=0.\qedhere
\edm
\end{proof}
\begin{NB}
For $s=y=1$, all four $3$-forms $\T_3,\ldots,\T_6$ coincide, reflecting
the fact that the undeformed metric has $\psi_3,\ldots,\psi_6$ as parallel
spinors for the connection defined by
\bdm
\T\ :=\ \frac{1}{4}\,[X_{135}\, +\, X_{146}\, +\,X_{245}\, -\,X_{236}]\,-\,
\frac{1}{2}\,X_{127}\, - \, \frac{1}{4}\left[ X_{347}-X_{567}\right]\,.
\edm
The $3$-forms $\T_3$ and $\T_4$ are equal for the family of metrics
defined by $2s^2=1+y$, whereas $\T_5=\T_6$ as soon as $y=1$. Even more
interestingly, there exists a metric for which $\T_3=-\T_4$:
\end{NB}
\begin{prop}\label{TandminusT}
Consider the metric $g_{s_0,y_0}$ on $N(1,1)$ defined by $s_0=\sqrt{3}/2$ and 
$y_0=2$, and the $3$-form
\bdm
\T \ :=\ \sqrt{3}/6\,(X_{135}+X_{146}-X_{245}+X_{236}).
\edm
Then, $\psi_3$ is parallel with respect to the connection $\nabla^{4 \cdot 
\T}$, and $\psi_4$ is parallel with respect to the connection $\nabla^{- 
4 \cdot \T}$.
Furthermore, both connections are not flat.
\end{prop}
\noindent
It is a subtle and computationally difficult question in as much
$\T$ can be adapted to a given spinor in order to make it parallel.
For this, a more systematic approach is required. There are precisely $13$
isotropy invariant $3$-forms on $\m$, hence the most general $3$-form
we can consider is a linear combination of
\bdm
X_{135}+X_{146},\ X_{235}+X_{246},\ X_{357}+X_{467},\ X_{145}-X_{136},\
X_{245}-X_{236},\ X_{457}-X_{367},
\edm\bdm
X_{127},\ X_{347}, \ X_{567},\ 
X_{134},\ X_{234},\ X_{156},\ X_{256}.
\edm
We studied the question whether there exists a continuous family of $3$-forms 
$\T_{a,b}$ of this general type such that a given linear combination
$a\,\psi_3+b\,\psi_5$ is parallel with respect to $\nabla^{\T_{a,b}}$.
It turns out that this is possible if and only if $s=y$.
We state the result of this lengthy calculation without proof.
\begin{prop} \label{Tab} 
The spinor field $\psi_{a,b}:=a\cdot \psi_3+b\cdot \psi_5,\, ab\neq 0$, 
satisfies
the equation $\nabla^g_{X}\psi_{a,b}+(X\haken\T_{a,b})\cdot\psi_{a,b} = 0$ 
if and only if $s=y$ and if $\T_{a,b}$ is given by
\begin{eqnarray*}
\T_{a,b} & = & \frac{a^2(-7s^2+8s+2)+ b^2(s^2+4s-2)}{12s(a^2+b^2)}\, 
[X_{135}+X_{146}] +\frac{s^2+4s-2}{12s}[X_{245}-X_{236}] \\
&+&\frac{a^2(-8s^2+s+4)+b^2(-4s^2+5s-4)}{12s(a^2+b^2)}[X_{347}-X_{567}]
+ \frac{-4s^2+2s-1}{6s}X_{127}\\
&+& \frac{ab\,(-2s^2+s+1)}{3s(a^2+b^2)}[X_{134}-X_{156}] +
\frac{ab\,(s^2+s-2)}{3s(a^2+b^2)}[X_{357}+X_{467}]
\end{eqnarray*}
\end{prop}
\noindent
For $ab=0$, this differential form $\T_{a,b}$ is again a linear combination
of the seven basic $3$-forms we started with, and coincides indeed
for $a=1,b=0$ and $a=0,b=1$ with the $3$-forms $\T_3,\,\T_5$ 
evaluated at the parameter value $s=y$, respectively. Remark that 
the connections with torsion $\T_{a,b}$ constitute a $S^1$-parameter family of
connections admitting parallel spinors on the same Riemannian manifold.  
The $3$-Sasakian metric ($s=1,y=2$) and the Einstein metric ($s=1,y=2/5$)
are of particular interest. For theoretical reasons to be explained
in the next section, both must admit a family of torsion forms such that
the three Killing spinors of the $3$-Sasakian metric ($\psi_3,\psi_4,\psi_6$
in our notation) are parallel with respect to the connection it defines.
In fact, such a family exists for $s=1$ and arbitrary $y$ (but not for
arbitrary $s$).
\begin{prop}
For the metrics $g_{s,y}$ on $N(1,1)$, the spinor field
$\psi_{a,b,c}:=a\,\psi_3+b\,\psi_4+c\,\psi_6,\, abc\neq0$, satisfies the equation
$\nabla^g_X\psi_{a,b,c}+ (X\haken \T_{a,b,c})\cdot\psi_{a,b,c}=0$
if and only if $s=1$ and if $\T_{a,b,c}$ is given by 
\begin{eqnarray*}
\T_{a,b,c}\!\!\!& =&\!\!\!  
P(a,b,c)[X_{567}- X_{347}]+ P(a,c,b)[X_{135}+X_{146}]
+ P(b,c,a)[X_{245}- X_{236}]\\
\!\!\!&+&\!\!\! Q(a,b,c)[X_{235}+X_{246}+X_{145}-X_{136}]
+ Q(b,c,a)[X_{357}+X_{467}+X_{156}-X_{134}]\\
\!\!\!&+& \!\!\!Q(a,c,b)[X_{457}-X_{367}+X_{234}-X_{256}] 
+ \frac{2y-5}{6}\,X_{127},
\end{eqnarray*}
with the following definitions for the coefficients $P$ and $Q$:
\bdm
P(a,b,c)\ :=\ \frac{(a^2+b^2)(4-y)+c^2(8-5y)}{12(a^2+b^2+c^2)},\quad
Q(a,b,c)\ :=\ \frac{ab\,(y-1)}{3(a^2+b^2+c^2)}.
\edm
\end{prop}

\noindent
Let us discuss the spinor fields $\psi_3$ and $\psi_5$ from the point of 
view of $\mathrm{G}_2$-geometry.
In general, a spinor field $\psi$ of length one defines on a $7$-dimensional 
Riemannian manifold a $3$-form of general type by the formula
(see \cite{BFGK}, \cite{FKMS})
\bdm
\omega(X,Y,Z) \ := \ - \, \langle X \cdot Y \cdot Z \cdot \psi \, , 
\, \psi \rangle \ . 
\edm 
Computing the forms of the spinors $\psi_3, \, \psi_5$ we obtain
\begin{eqnarray*}
\omega_3 & = & - \, X_{127} \, + \, X_{135} \, + X_{146} \, + \, X_{236}
\, - \, X_{245} \, - \, X_{347} \, + \, X_{567} \ , \\ 
\omega_5 & = & - \, X_{127} \, - \, X_{135} \, - X_{146} \, + \, X_{236}
\, - \, X_{245} \, + \, X_{347} \, - \, X_{567} \ .
\end{eqnarray*} 
The connections $\nabla^{3}$ and $\nabla^{5}$ with torsion forms 
$4 \cdot \T_3$ and $4 \cdot \T_5$ 
preserve the 
$\mathrm{G}_2$-structures $\omega_3$ and $\omega_5$, respectively.
Moreover, a direct computation yields the formulas
\begin{eqnarray*} 
\big( \T_3 \, , \, \omega_3 \big) & = & \frac{4 \cdot s^2 + y + 1}
{6 \cdot s} \ > \ 0 \, ,
\\
\big( \T_5 \, , \, \omega_5 \big) & = & - \, \frac{4 \cdot s + 2 \cdot s \cdot y + y - 1}{6 \cdot s} \ .
\end{eqnarray*} 
Since the connection preserving a $\mathrm{G}_2$-structure
is unique (see \cite{Friedrich&I1}), the $\mathrm{G}_2$-structures
$\omega_3$ and $\omega_5$ are not equivalent. We remark that
$\omega_3$ and $\omega_5$ are cocalibrated $\mathrm{G}_2$-structures,
\bdm
d * \omega_3 \ = \ 0 \, , \quad d * \T_3 \ = \ 0 \, , \quad 
d * \omega_5 \ = \ 0 \, , \quad d * \T_5 \ = \ 0 \ .
\edm 
Indeed, for any vector, the inner product
$X \haken * \, \omega_3$ is orthogonal to $7 \cdot \T_3 - (\T_3 \, , \, 
\omega_3) 
\cdot \omega_3$. The formula expressing the torsion form $\T$ of an admissible
$\mathrm{G}_2$-structure by the $3$-form $\omega$  
( see \cite{Friedrich&I1} and \cite{Friedrich&I3}) yields now 
$d * \omega_3 = 0$ immediately. 
The codifferential of the torsion form is given by the formula
(see \cite{Friedrich&I1} and \cite{Friedrich&I3})
\bdm
4 \cdot d * \T_3 \ = \ d \lambda_3 \wedge * \omega_3 \, , \quad
\lambda_3 \ = \ (4 \cdot \T_3 \, , \, \omega_3) \ .
\edm
In our example the function $\lambda_3$ is constant, i.e., $d*\T_3 = 0$.
The same argument
applies for $\omega_5$. The class of all cocalibrated $\mathrm{G}_2$-structures splits into the sum $\mathcal{W}_1 \oplus \mathcal{W}_3$ of a 
$1$-dimensional class $\mathcal{W}_1$ (the so called nearly parallel
$\mathrm{G}_2$-structures) and a $27$-dimensional class $\mathcal{W}_3$
(see \cite{FG}). Nearly parallel $\mathrm{G}_2$-structures are characterized 
by the
condition that the torsion form $\T$ of its unique connection is proportional
to $\omega$. On the other side, the 
$\mathrm{G}_2$-structures of type $\mathcal{W}_3$ are the cocalibrated 
structures such that $\T$ and $\omega$ are orthogonal, $(\T \, , \, \omega) 
= 0$ (see \cite{Friedrich&I1}). Using this
characterization we obtain immediately
\begin{prop}
The $\mathrm{G}_2$-structure $\omega_3$ is nearly parallel if and
only if $s = 1$ and $y = 2$.
The $\mathrm{G}_2$-structure $\omega_3$ is never of type $\mathcal{W}_3$.
The $\mathrm{G}_2$-structure $\omega_5$ is nearly parallel if and
only if $s = 1$ and $y = 2/5$. This metric is a universal deformation of the
$3$-Sasakian metric (see \cite{FKMS}).
The $\mathrm{G}_2$-structure $\omega_5$ is of type $\mathcal{W}_3$
if and only if $ 2 \cdot s \cdot (2 + y) = 1 - y$.
\end{prop} 
\noindent 
In general, the scalar curvatures $\mathrm{Scal}^g \, , \, \mathrm{Scal}^{\nabla}$  of a cocalibrated $\mathrm{G}_2$-structure $(M^7, g, \omega)$ can be expressed by its torsion form $\T$ 
(see \cite{Friedrich&I3}) :
\bdm
\mathrm{Scal}^g \ = \ 2 \cdot (\T \, , \, \omega)^2 \, - \, \frac{1}{2} \cdot
||\T||^2 \, , \quad \mathrm{Scal}^{\nabla} \ = \ \mathrm{Scal}^g \, - \, \frac{3}{2} \cdot ||\T||^2 \ = \ 2 \cdot (\T \, , \, \omega)^2 \, - \, 2 \cdot
||\T||^2 \ .
\edm
We use the forms $\omega_3 \, , \,  4 \cdot \T_3$ as well as the forms
$\omega_5 \, , \, 4 \cdot \T_5$ in order to compute the Riemannian scalar curvature of the metric depending on the parameters $s,y$. In both cases the result
is the same :
\bdm
\mathrm{Scal}^g \ = \ 8 \, + \, 24 y \, - \, 2 y^2 \, - \, 
\frac{2 + y^2}{s^2} \ .
\edm
In a similar way we compute the scalar curvature of the connection
$\nabla^{3}$ and $\nabla^{5}$ :
\bdm
\mathrm{Scal}^3 \ = \ - \, \frac{4}{3 s^2}\cdot \big(8 + 32 s^4 + 4 y + 
5 y^2 + 2 s^2 (-4 - 28 y + 3 y^2)\big) \ ,
\edm
\bdm
\mathrm{Scal}^5 \ = \ - \, \frac{4}{3 s^2} \cdot \big(8 - 4 y + 5 y^2 + 
8 s (-2 + y + y^2) + 2 s^2 (4 - 20 y + 7 y^2)\big) \ .
\edm
In particular, we obtain a family of cocalibrated $\mathrm{G}_2$-structures
on $N(1,1)$ with vanishing scalar curvature of the associated connection. 
Moreover, a numerical computation yields that there exist two pairs
of parameters where both scalar curvatures $\mathrm{Scal}^3$ and $\mathrm{Scal}^5$ vanish, namely $(s\, , \, y) \approx ( 0.62066 \, , \,  0.852508)$ and 
$(1.49934 \, , \, 1.66564)$. The Ricci tensor $\mathrm{Ric}^{\nabla}$ of 
the canonical connection of a $\mathrm{G}_2$-structure $(M^7, g , \omega)$
can be expressed by the derivative $d \T$ of the torsion form (see \cite{Friedrich&I1}),
\bdm
\mathrm{Ric}^{\nabla}(X_i \, , \, X_j) \ = \ \frac{1}{2} \cdot \big (X_i
\haken d \T \, , \, X_j \haken * \omega \big) \ .
\edm
Using the commutator relations in the Lie algebra we compute the
exterior derivatives 
\begin{eqnarray*}
d X_1 &=& - \, 2 s \cdot X_{27}  \, + \, y \cdot (X_{35} \, + \, X_{46}) \, ,\\
d X_2 &=& \ \ 2 s \cdot X_{17}  \, + \, y \cdot (X_{45} \, - \, X_{36}) \, ,\\
d X_7 &=& - \, \frac{2}{s} \cdot X_{12}  \, - \, \frac{y}{s} \cdot (X_{34} \, - \, X_{56}) \ .
\end{eqnarray*} 
The torsion form $\T_3$ can be written as
\begin{eqnarray*}
\T_3 &=& \ - \, \frac{(y-2)(5s^2 - 1 - y)}{3 s y}  X_{127} \, + \,
\frac{1}{y} \Big[ \frac{1}{2} \, - \, \frac{y}{4} \, + \, \frac{1 + y}{6 s} 
\, - \, \frac{s}{3} \Big]  X_1 \wedge d X_1 \\
 &+&  \frac{1}{y} \Big[ \frac{1}{2} \, - \, \frac{y}{4} \, - \, 
\frac{1 + y}{6 s} \, + \, \frac{s}{3} \Big] X_2 \wedge d X_2 \, 
- \, \frac{s}{y} \Big[\frac{4 + y}{12 s} \, - \, \frac{2 s}{3}\Big] 
X_7 \wedge d X_7 \ .
\end{eqnarray*}
We can now compute the exterior derivative as well as the Ricci tensor.
Let us discuss the cases $(s,y) = (1,4)$ and $(s,y) = 
(\sqrt{3/2}\, , \, 2)$
where the formulas simplify.
\begin{exa}\label{scal-null} 
In case of $s = 1$ and $y = 4$ we obtain
\bdm
\T_3 \ = \ - \, \frac{1}{4} X_2 \wedge d X_2 \, ,
\quad d \T_3 \ = \ 8  X_{3456} \, - \, 4  X_{1457} \, + \, 
4  X_{1367}  .
\edm
The scalar curvatures are  $\mathrm{Scal}^3 = 0$ and 
$\mathrm{Scal}^g = 54$ . Moreover, we obtain
\bdm
- \, 2 \cdot \sigma_{\T_3} \ = \ (\T_3^2)_4 \  =  \ \frac{1}{4} 
\cdot d \T_3 ,
\edm
i.e., $d\T_3$ is proportional to $\sigma_{\T_3}$ (see Theorem \ref{parameter}).
\end{exa}
\begin{exa} In case of $s = \sqrt{3/2}$ and $y = 2$ we obtain
\bdm
\T_3 \ = \ \frac{1}{4} X_7 \wedge d X_7 \, ,
\quad d \T_3 \ = \ \frac{4}{3} (X_{1234} \, - \, X_{1256} \, - \, 
X_{3456}) .
\edm
The scalar curvatures $\mathrm{Scal}^3$ and 
$\mathrm{Scal}^g$ are positive.
\end{exa}
\vspace{2mm}

\noindent
The $3$-form $\omega_{a,b}$ corresponding to the 
spinor $a \cdot \psi_3 + b \cdot \psi_5$ is given by the formula
\begin{eqnarray*}
\omega_{a,b} \ = &-& (a^2 + b^2)( X_{127} \, - \, X_{236} \, 
+ \, X_{245} ) \, 
+ \, 2 a b 
(X_{134} \, - \, X_{156} \, + \, X_{357} \, + \, X_{467}) \\
&+& (a^2 - b^2) (X_{135} \, + \, X_{146} \, - \, X_{347} \, + \, X_{567})  . 
\end{eqnarray*} 
We compute the inner product with the torsion form
of Proposition \ref{Tab} :
\bdm
( \T_{a,b} \, , \, \omega_{a,b}) \ = \ \frac{b^2 (1 - 5 \, s - 2 
\, s^2) + a^2 (1 + s + 4 \, s^2)}{6 \, s}\,  .
\edm
In particular, the $\mathrm{G}_2$-structure is of pure type $\mathcal{W}_3$
if and only if
\bdm 
b^2 (1 - 5 \, s - 2 \, s^2) + a^2  (1 + s + 4 \, s^2) \ = \ 0 \, , \quad 
y \ = \ s \, .
\edm
\noindent
Finally, we will construct non trivial $4$-forms on $N(1,1)$ such that the 
underlying connections admit parallel spinors. Remark that
spinorial connections related to $4$-forms are of completely different
type. For example, they do not preserve the hermitian product of 
spinors and, in general, the holonomy group of a spinorial connection of that
type is non compact. Nevertheless, for the family of metrics $g_{s,y}$,
the qualitative results are quite similar to those for $3$-forms, though
they cannot be deduced one from each other.
\begin{thm}
On $N(1,1)$ with the metric $g_{s,y}$, there exist four
spinorial connections defined by $4$-forms admitting $\psi_3,\psi_4,\psi_5$
 and 
$\psi_6$ as parallel spinors, respectively;  and there exists, for $s=1$,
arbitrary $y$ and any linear combination of $\psi_3,\psi_4,\psi_6$
precisely one $4$-form which makes this particular combination parallel.
\end{thm}
\noindent
Again, the  exposition of results will make the statement more precise. 
No proofs will be given, since they are similar to the corresponding
computations for $3$-forms. 
With the same notations as before, consider now the following Ansatz for a 
global $4$-form :
\begin{eqnarray*}
R& =& \alpha\, X_{1234}+\beta\,X_{1256}+\gamma\,X_{3456}+\delta\,X_{1347}+
\eps\,X_{1567}+\xi\,X_{2347}+\mu\,X_{2567}\\
&+&\nu\,(X_{1235}+X_{1246})+ \lambda\,(X_{1357}+X_{1467})
+\eta\,(X_{1245}-X_{1236}) \\ 
&+&\omega\,(X_{1457}-X_{1367})
+ \pi\,(X_{2457}-X_{2367})+\vrho\,(X_{2357}+X_{2467}).
\end{eqnarray*}
\begin{prop}
The spinor field $\psi_3$ satisfies the equation 
$\nabla^g_X \psi_3+ (X\haken R)\cdot \psi_3=0$ if and only if $R=R_3$, with
\begin{eqnarray*}
R_3 & := &
\left[\frac{y}{4}-\frac{1}{2}+\frac{1+y}{8s}-\frac{s}{2}\right][X_{1457}-X_{1367}]
+\left[\frac{1}{2}-\frac{y}{4}+\frac{1+y}{8s}-\frac{s}{2}\right] [X_{2357}
+X_{2467}]\\
&+& \left[-\frac{s}{2}+\frac{y+3}{8s} \right][X_{1234} -X_{1256}] +
\left[\frac{s}{2}+\frac{1-3y}{8s}\right] X_{3456}.
\end{eqnarray*}
\end{prop}
\begin{prop}
The spinor field $\psi_4$ satisfies the equation 
$\nabla^g_X \psi_4+ (X\haken R)\cdot \psi_4=0$ if and only if $R=R_4$, with
\begin{eqnarray*}
R_4 & := &
\left[\frac{y}{4}-\frac{1}{2}-\frac{1+y}{8s}+\frac{s}{2}\right][X_{1457}-X_{1367}]
+\left[\frac{1}{2}-\frac{y}{4}-\frac{1+y}{8s}+\frac{s}{2}\right] [X_{2357}
+X_{2467}]\\
&+& \left[-\frac{s}{2}+\frac{y+3}{8s} \right][X_{1234} -X_{1256}] +
\left[\frac{s}{2}+\frac{1-3y}{8s}\right] X_{3456}.
\end{eqnarray*}
\end{prop}
\begin{prop}
The spinor field $\psi_5$ satisfies the equation 
$\nabla^g_X \psi_5+ (X\haken R)\cdot \psi_5=0$ if and only if $R=R_5$, with
\begin{eqnarray*}
R_5 & := &
 \left[\frac{1}{2}-\frac{y}{4}+\frac{y-3}{8s} \right][X_{1234} -X_{1256}] +
\left[-\frac{1}{2}+\frac{3y}{4}+\frac{1+3y}{8s}\right] X_{3456}\\
&+&\frac{1-y}{8s}[(X_{1457}-X_{1367})-(X_{2357} +X_{2467})].
\end{eqnarray*}
\end{prop}
\begin{prop}
The spinor field $\psi_6$ satisfies the equation 
$\nabla^g_X \psi_6+ (X\haken R)\cdot \psi_6=0$ if and only if $R=R_6$, with
\begin{eqnarray*}
R_6 & := &
 \left[-\frac{1}{2}+\frac{y}{4}+\frac{y-3}{8s} \right][X_{1234} -X_{1256}] +
\left[\frac{1}{2}-\frac{3y}{4}+\frac{1+3y}{8s}\right] X_{3456}\\
&+&\frac{y-1}{8s}[(X_{1457}-X_{1367})-(X_{2357} +X_{2467})].
\end{eqnarray*}
\end{prop}
\begin{NB}
For $s=y=1$, no two of the  four $4$-forms $R_3,\ldots,R_6$ coincide, 
reflecting the different behavior of spinorial connections defined by
$4$-forms when compared to connections defined by $3$-forms.
The $4$-forms $R_3$ and $R_4$ are equal for the family of metrics
defined by $4s^2=1+y$, whereas $R_5$ and $R_6$ are never equal. 
As for $3$-forms, there exists a metric for which  $R_3=-R_4$:
\end{NB}
\begin{prop}\label{RandminusR}
Consider the metric $g_{s_0,y_0}$ on $N(1,1)$ defined by $s_0=\sqrt{5}/2$ and 
$y_0=2$, and the $4$-form
\bdm
R \ :=\ -\sqrt{5}/10\,[(X_{1457}-X_{1367})+(X_{2357}+X_{2467})].
\edm
Then, $\psi_3$ is parallel with respect to the connection $\nabla^{R}$,
and $\psi_4$ is parallel with respect to the connection $\nabla^{-R}$.
Furthermore, both connections are not flat.
\end{prop}
\begin{prop}
For the metrics $g_{s,y}$ on $N(1,1)$, the spinor field
$\psi_{a,b,c}:=a\,\psi_3+b\,\psi_4+c\,\psi_6,\, abc\neq0$, satisfies the equation
$\nabla^g_X\psi_{a,b,c}+ (X\haken R_{a,b,c})\cdot\psi_{a,b,c}=0$
if and only if $s=1$ and if $R_{a,b,c}$ is given by 
\begin{eqnarray*}
R_{a,b,c}\!\!\!& =&\!\!\!  
P(a,b,c)[X_{1234}- X_{1256}]- P(a,c,b)[X_{2467}+X_{2357}]
+ P(b,c,a)[X_{1457}- X_{1367}]\\
\!\!\!&+&\!\!\! Q(a,b,c)[(X_{2457}-X_{2367})-(X_{1357}+X_{1467})]
+ Q(b,c,a)[(X_{1235}+X_{1246})+\\
\!\!\!&+& \!\!\! (X_{2567}-X_{2347})]
 + Q(a,c,b)[(X_{1567}-X_{1347})-(X_{1245}-X_{1236})] 
+ \frac{5-3y}{8}\,X_{3456},
\end{eqnarray*}
with the following definitions for the coefficients $P$ and $Q$:
\bdm
P(a,b,c)\ :=\ \frac{(a^2+b^2)(y-1)+c^2(3y-7)}{8(a^2+b^2+c^2)},\quad
Q(a,b,c)\ :=\ \frac{ab\,(y-3)}{4(a^2+b^2+c^2)}.
\edm
\end{prop}
%
\section{Torsion forms with parallel spinors on $3$-Sasakian manifolds}
\noindent
The Aloff-Wallach space $N(1,1)$ admits a $3$-Sasakian structure, and
some special torsion forms with parallel spinors discussed in Section $8$
are closely related to the underlying contact structures of $N(1,1)$. This 
observation yielded the idea that \emph{any} $3$-Sasakian manifold 
should admit natural 
connections with skew-symmetric torsion and parallel spinors. In this 
section, we will make this remark precise. In particular, for
a fixed $3$-Sasakian metric, we will construct a whole family of connections 
with parallel spinors. The structure group of a $3$-Sasakian geometry is 
the subgroup $\SU(2) \subset \mathrm{G}_2 \subset \SO(7)$,
the isotropy group of four spinors in dimension seven. In order to
keep the realization of the spin representation we used in Section $8$,
we describe the subgroup $\SU(2)$ in such a way that the vectors
$e_1, e_2, e_7 \in \R^7$ are fixed. More precisely, the Lie algebra
$\su(2)$ is generated by the following $2$-forms in $\R^7$:
\bdm
e_{34} \, + \, e_{56} \, , \quad 
e_{35} \, - \, e_{46} \, , \quad
e_{36} \, + \, e_{56} \ .
\edm
The real spin representation $\Delta_7$
splits under the action of $\SU(2)$ into a $4$-dimensional trivial
representation $\Delta_7^0$ and the unique non trivial $4$-dimensional
representation $\Delta_7^{1}$. In our spin basis, the space $\Delta_7^0$
is spanned by the spinors $\psi_3, \psi_4, \psi_5, \psi_6$. We consider
the following $\SU(2)$-invariant $2$-forms on $\R^7$:
\bdm
de_1 \ := \ e_{35} \, + \, e_{46} \, , \quad 
de_2 \ := \ e_{45} \, - \, e_{36} \, , \quad
de_7 \ := \ e_{34} \, - \, e_{56} \ .
\edm
Using this notation, we introduce a family of invariant $3$-forms in 
$\R^7$ depending on $10$ parameters,
\bdm
\T \ = \ \sum_{i,j = 1,2,7} x_{ij} \cdot e_i \wedge d e_j \, + \, w \cdot e_1 \wedge e_2 
\wedge e_7
\edm
The key point of our considerations in this section is the following 
algebraic observation
\begin{prop}\label{3Sas}  
For any spinor $\psi \in \Delta_7^0$, there exists a unique
invariant $3$-form $\T$ such that 
$\big\{ X -  2 \cdot X \haken \T \big\} \cdot \psi = 0$ holds for any 
vector $X \in \R^7$. 
\end{prop} 
\begin{proof} Given a spinor $\psi = a \, \psi_3 + b \, \psi_4 + 
c \, \psi_5 + d \, \psi_6$, we solve the overdetermined system
$(X - 2 \cdot X \haken \T) \cdot \psi = 0$ with respect to the coefficients
of the $3$-form. It turns out that a solution exists and is given by
the following explicit formulas
\begin{eqnarray*}
x_{11} & = & \frac{a^2-b^2-c^2+d^2}{6(a^2+b^2+c^2+d^2)} , \ 
x_{12} \ = \ \frac{ab+cd}{3(a^2+b^2+c^2+d^2)} , \
x_{17} \ = \ \frac{ac - bd}{3(a^2+b^2+c^2+d^2)} , \\
x_{21} & = & \frac{ab - cd}{3(a^2+b^2+c^2+d^2)} , \ 
x_{22} \ = \ \frac{-a^2+b^2-c^2+d^2}{6(a^2+b^2+c^2+d^2)} , \ 
x_{27} \ = \ \frac{bc + ad}{3(a^2+b^2+c^2+d^2)} , \\
x_{71} & = & \frac{ac + bd}{3(a^2+b^2+c^2+d^2)} , \
x_{72} \ = \ \frac{bc-ad}{3(a^2+b^2+c^2+d^2)} , \ 
x_{77} \ = \ \frac{-a^2-b^2+c^2+d^2}{6(a^2+b^2+c^2+d^2)} .
\end{eqnarray*}
and $ w = - \, 1/6$. The map $(a,b,c,d) \rightarrow x_{ij}(a,b,c,d)$ is
the Veronese map from $\mathbb{P}^3$ into the sphere $S^8$ of radius $1/12$. 
\end{proof}
\noindent
Consider a simply connected $3$-Sasakian manifold $M^7$ of dimension seven and 
denote its three contact structures by $\eta_1, \eta_2$, and $\eta_7$. 
It is known that $M^7$ is then an Einstein space,  and
examples (also non homogeneous ones) can be found in the paper \cite{BG}
by Boyer and Galicki. The tangent bundle of $M^7$ splits
into the $3$-dimensional part spanned by $\eta_1, \eta_2, \eta_7$ and
its $4$-dimensional orthogonal complement. We restrict the exterior derivatives $d \eta_1, d \eta_2$ and $d \eta_7$ to this complement. In an adapted
orthonormal frame, these forms coincide with the algebraic forms $de_1, de_2$
and $de_7$. Now we apply Proposition \ref{3Sas}. The space 
of Riemannian Killing spinors 
\bdm
\nabla^g_X \psi \, + \, \frac{1}{2} \cdot X \cdot \psi \ = \ 0
\edm
is non trivial and has at least dimension three (see \cite{FK}). Moreover,
the proof of this fact shows that all the Riemannian Killing spinors
are sections in the subbundle corresponding to the $\SU(2)$-representation
$\Delta_7^0$. Consequently, for any Killing spinor, there exists a unique
torsion form $\T$ of the described type such that
\bdm
\nabla_X^{\T} \psi \ = \ \nabla^g_X \psi \, + \, (X \haken \T) \cdot \psi 
\ = \ 0 \ .
\edm
\begin{thm} Any $3$-Sasakian manifold in dimension seven admits a
$\mathbb{P}^2$-parameter family of metric connections with
skew-symmetric torsion and parallel spinors. The holonomy group of these
connections is a subgroup of $\mathrm{G}_2$.
\end{thm}
\noindent
The space of $\SU(2)$-invariant $4$-forms on $\R^7$ has dimension ten,
\bdm
\T \ = \ \sum_{i,j,k=1,2,7} x_{ijk} \cdot e_i \wedge e_j 
\wedge d e_k \, + \, w \cdot e_3 \wedge e_4 \wedge e_5 \wedge e_6 \ .
\edm
We study spinorial connections depending on $4$-forms. Again, any spinor
in $\Delta_7^0$ defines a unique $4$-form being a solution of the
corresponding overdetermined linear system and we can apply the
same construction as above. Let us formulate the results. 
\begin{prop} For any spinor $\psi \in \Delta_7^0$ there exists a unique
invariant $4$-form $\T$ such that $\big\{ X - 2 \cdot X \haken \T \big\} \cdot \psi = 0$ holds for any vector $X \in \R^7$.
\end{prop} 
\begin{thm} \label{3Sas4form}  
Any $3$-Sasakian manifold in dimension seven admits a
$\mathbb{P}^2$-parameter family of spinorial connections defined by $4$-forms 
and with parallel spinors. The spinorial holonomy group of these connections
is a subgroup of $\GL(7, \R)$.
\end{thm}

\appendix
\section{The Lie algebra $\spin(9)$ inside $\so(16)$}
%
\noindent
The Lie algebra $\so(16)$ of all antisymmetric matrices is parameterized
by $120$ parameters $\omega_{i,j}$, $1 \leq i < j \leq 16$. We realize
the $36$-dimensional subalgebra $\spin(9)$ by $84$ explicit equations. 
The first group of $56$ equations involves forms of type 
$\omega_{8,\alpha}$ and 
$\omega_{i,\beta}$, where $1 \leq i < 8 < \alpha, \beta \leq 16$,
and is given in Table~\ref{table-1}.
\begin{table}
\bdm\begin{array}{|l|l|l|l|l|}
\hline
\omega_{1,9} =  \omega_{8,16} & \omega_{1,10} = -\omega_{8,15} &
\omega_{1,11} =  \omega_{8,14}& \omega_{1,12} = \omega_{8,13} &
\omega_{1,13} = -\omega_{8,12} \\[1mm] \hline
\omega_{1,14} = - \omega_{8,11} & \omega_{1,15} = \omega_{8,10} & 
\omega_{1,16} = -\omega_{8,9}  & \omega_{2,9} =  \omega_{8,15} & 
\omega_{2,10} = \omega_{8,16} \\[1mm] \hline
\omega_{2,11} = -\omega_{8,13} & \omega_{2,12} = \omega_{8,14} &
\omega_{2,13} = \omega_{8,11} & \omega_{2,14} = - \omega_{8,12} &
\omega_{2,15} = -\omega_{8,9} \\ [1mm] \hline
\omega_{2,16} = -\omega_{8,10}& \omega_{3,9} = - \omega_{8,14} &
\omega_{3,10} =  \omega_{8,13}& \omega_{3,11} =  \omega_{8,16} &
\omega_{3,12} =  \omega_{8,15} \\ [1mm] \hline
\omega_{3,13} = -\omega_{8,10} & \omega_{3,14} = \omega_{8,9} &
\omega_{3,15} = -\omega_{8,12} & \omega_{3,16} = -\omega_{8,11} &
\omega_{4,9}  = -\omega_{8,13} \\ [1mm] \hline
\omega_{4,10} = -\omega_{8,14} & \omega_{4,11} = -\omega_{8,15} &
\omega_{4,12} =  \omega_{8,16} & \omega_{4,13} =  \omega_{8,9} &
\omega_{4,14} = \omega_{8,10} \\ [1mm] \hline
\omega_{4,15} = \omega_{8,11} & \omega_{4,16} = - \omega_{8,12} &
\omega_{5,9}  = \omega_{8,12} & \omega_{5,10} = - \omega_{8,11} &
\omega_{5,11} = \omega_{8,10} \\ [1mm] \hline
\omega_{5,12} = -\omega_{8,9} & \omega_{5,13} =  \omega_{8,16} &
\omega_{5,14} = \omega_{8,15} & \omega_{5,15} = -\omega_{8,14} &
\omega_{5,16} = - \omega_{8,13} \\ [1mm] \hline
\omega_{6,9}  = \omega_{8,11} & \omega_{6,10} =  \omega_{8,12} &
\omega_{6,11} = -\omega_{8,9} & \omega_{6,12} = - \omega_{8,10} &
\omega_{6,13} = -\omega_{8,15}\\ [1mm] \hline
\omega_{6,14} =  \omega_{8,16} & \omega_{6,15} = \omega_{8,13} &
\omega_{6,16} = - \omega_{8,14}& \omega_{7,9}  = -\omega_{8,10} &
\omega_{7,10} =   \omega_{8,9} \\ [1mm] \hline
\omega_{7,11} =  \omega_{8,12} & \omega_{7,12} = - \omega_{8,11} & 
\omega_{7,13} =  \omega_{8,14} & \omega_{7,14} = - \omega_{8,13} &
\omega_{7,15} =  \omega_{8,16} \\ [1mm] \hline
\omega_{7,16} = -\omega_{8,15} & & & & \\ [1mm] \hline
\end{array}\edm
\caption{First group of equations defining $\spin(9)$ inside $\so(16)$.}
\label{table-1}
\end{table}
The second group of $28$ equations involves the forms $\omega_{i,j} \, , \, 
\omega_{\alpha, \beta}$ for $1 \leq i,j \leq 8 <  \alpha, \beta 
\leq 16$, and is given in Table~\ref{table-2}. 
\begin{table}
\bdm
\begin{array}{|l|l|}\hline
2 \cdot \omega_{1,2} = \omega_{11,12} +  \omega_{13,14} -  
\omega_{15,16} +  \omega_{9,10} &
2 \cdot \omega_{1,3} =  -  \omega_{10,12} +  \omega_{13,15} +  
\omega_{14,16} +  \omega_{9,11} \\[2mm] \hline
2 \cdot \omega_{1,4} =  \omega_{10,11} +  \omega_{13,16} -  
\omega_{14,15} +  \omega_{9,12} &
2 \cdot \omega_{1,5}  =  -  \omega_{10,14} -  \omega_{11,15} -  
\omega_{12,16} +  \omega_{9,13} \\[2mm] \hline
2 \cdot \omega_{1,6}  =  \omega_{10,13} -  \omega_{11,16} + 
\omega_{12,15} +  \omega_{9,14} &
2 \cdot \omega_{1,7} \ = \ \omega_{10,16} + \omega_{11,13} -  
\omega_{12,14} + \omega_{9,15}  \\[2mm] \hline
2 \cdot \omega_{1,8}  =  - \omega_{10,15} +  \omega_{11,14} + 
\omega_{12,13} +  \omega_{9,16} &
2 \cdot \omega_{2,3} \ = \ \omega_{10,11} - \omega_{13,16} +  
\omega_{14,15} + \omega_{9,12}  \\[2mm] \hline
2 \cdot \omega_{2,4}  =  \omega_{10,12} +  \omega_{13,15} + 
\omega_{14,16} -  \omega_{9,11} &
2 \cdot \omega_{2,5} \ = \ \omega_{10,13} + \omega_{11,16} -  
\omega_{12,15} + \omega_{9,14}  \\[2mm] \hline
2 \cdot \omega_{2,6}  =  \omega_{10,14} -  \omega_{11,15} - 
\omega_{12,16} -  \omega_{9,13} &
2 \cdot \omega_{2,7} \ = \ \omega_{10,15} + \omega_{11,14} +  
\omega_{12,13} - \omega_{9,16}  \\[2mm] \hline
2 \cdot \omega_{2,8}  =  \omega_{10,16} -  \omega_{11,13} + 
\omega_{12,14} +  \omega_{9,15} &
2 \cdot \omega_{3,4} \ = \ \omega_{11,12} - \omega_{13,14} +  
\omega_{15,16} + \omega_{9,10}  \\[2mm] \hline
2 \cdot \omega_{3,5}  =  - \omega_{10,16} +  \omega_{11,13} + 
\omega_{12,14} +  \omega_{9,15} &
2 \cdot \omega_{3,6} \ = \ \omega_{10,15} + \omega_{11,14} -  
\omega_{12,13} + \omega_{9,16}  \\[2mm] \hline
2 \cdot \omega_{3,7}  =  - \omega_{10,14} +  \omega_{11,15} - 
\omega_{12,16} -  \omega_{9,13} &
2 \cdot \omega_{3,8} \ = \ \omega_{10,13} + \omega_{11,16} +  
\omega_{12,15} - \omega_{9,14}  \\[2mm] \hline
2 \cdot \omega_{4,5}  =  \omega_{10,15} -  \omega_{11,14} + 
\omega_{12,13} +  \omega_{9,16} &
2 \cdot \omega_{4,6} \ = \ \omega_{10,16} + \omega_{11,13} +  
\omega_{12,14} - \omega_{9,15}  \\[2mm] \hline
2 \cdot \omega_{4,7}  =  - \omega_{10,13} +  \omega_{11,16} + 
\omega_{12,15} +  \omega_{9,14} &
2 \cdot \omega_{4,8} \ = \ - \omega_{10,14} - \omega_{11,15} +  
\omega_{12,16} - \omega_{9,13}  \\[2mm] \hline
2 \cdot \omega_{5,6}  =  - \omega_{11,12} +  \omega_{13,14} + 
\omega_{15,16} +  \omega_{9,10} &
2 \cdot \omega_{5,7} \ = \ \omega_{10,12} + \omega_{13,15} -  
\omega_{14,16} + \omega_{9,11}  \\[2mm] \hline
2 \cdot \omega_{5,8}  =  - \omega_{10,11} +  \omega_{13,16} + 
\omega_{14,15} +  \omega_{9,12} &
2 \cdot \omega_{6,7} \ = \ \omega_{10,11} + \omega_{13,16} +  
\omega_{14,15} - \omega_{9,12}  \\[2mm] \hline
2 \cdot \omega_{6,8}  =  \omega_{10,12} -  \omega_{13,15} + 
\omega_{14,16} +  \omega_{9,11} &
2 \cdot \omega_{7,8} \ = \ \omega_{11,12} + \omega_{13,14} +  
\omega_{15,16} - \omega_{9,10}  \\[2mm] \hline
\end{array}
\edm
\caption{Second group of equations defining $\spin(9)$ inside $\so(16)$.}
\label{table-2}
\end{table}
\noindent
Consider a $3$-form $\T \in \T(\spin(9),\R^{16})$ in the
antisymmetric prolongation of the $\spin(9)$-representation in $\R^{16}$.
Then the $2$-forms $e_1 \haken \T  , \, e_8 \haken \T , \, 
e_9 \haken \T , \, e_{16} \haken \T$ are elements of $\spin(9)$. Using
the first equation $\omega_{1,9} = \omega_{8,16}$ defining this 
subalgebra, we conclude that
\bdm
\T_{1,8,9} \ = \ 0 \, , \quad \T_{1,8,16} \ = \ 0 \, , \quad 
\T_{8,9,16} \ = \ 0 \, , \quad \T_{1,9,16} \ = \ 0  .
\edm 
In a similar way, the first $56$ equations defining $\spin(9)$ yield
that, for $1 \leq i,j < 8$ and $8 < \alpha, \beta \leq 16$, the following
components of $\T$ vanish,
\bdm
\T_{i,8,\alpha} \ = \ 0 \, , \quad \T_{8,\alpha,\beta} \ = \ 0 .
\edm
The second $28$ equations immediately imply now that 
$\T_{i,j,8} = 0$,
i.e., the interior product $e_8 \haken \T = 0$ vanishes for {\it any}
$3$-form in the antisymmetric prolongation. Since the group $\mathrm{Spin}(9)$
acts transitively on the sphere in $\R^{16}$, we conclude that $\T = 0$.
\begin{prop} 
The antisymmetric prolongation of the unique irreducible $16$-dimensional 
of the Lie algebra $\spin(9)$  
vanishes, 
\bdm
\T(\spin(9), \R^{16}) \ = \ 0  .
\edm
\end{prop}
%
    
\end{document}